\documentclass[10pt]{article}

\usepackage{definitions}

\title{H-PDLP: GPU-Accelerated Linear Programming via Generalized PDHG and Homogeneous Self-Dual Embedding}
\author{Haihao Lu\thanks{MIT, Sloan School of Management (haihao@mit.edu).} \and Bo Tang\thanks{MIT, Sloan School of Management (botang@mit.edu).}}
\date{}

\begin{document}
\maketitle

\begin{abstract}
GPU-based first-order methods are emerging as an approach for solving large-scale linear programming (LP) problems. We present \HPDLP{}, a GPU first-order LP solver based on a generalized primal-dual hybrid gradient (PDHG) method for a homogeneous self-dual embedding of the LP, and study its numerical performance. The solver applies the restarted reflected Halpern iteration to the generalized PDHG operator, with the embedding encoding both primal-dual optimal solutions and certificates of primal or dual infeasibility. We derive a compact embedding for LPs with two-sided constraint and variable bounds in which the bounds enter through box cones, and we show that the resulting box-cone projections reduce to one-dimensional root-finding problems, so that each iteration remains matrix-free and well suited to GPUs. On the theoretical side, we establish conditions under which the iterates remain uniformly bounded away from the trivial zero solution of the embedding and derive complexity guarantees for the method. We also introduce \MIPLIBInfeas{}, a reproducible benchmark of 229 infeasible and 105 unbounded LPs derived from \MIPLIB{} 2017. Numerical experiments demonstrate the effectiveness of the proposed approach on both feasible LPs and infeasible or unbounded LPs.  The source code for \HPDLP{} is publicly available at \url{https://github.com/MIT-Lu-Lab/H-PDLP}.
\end{abstract}

\section{Introduction}
\label{sec:intro}

Linear programming (LP) is a foundational optimization model. Beyond its direct use in applications such as scheduling, network flow, chip design, and resource allocation, LP provides relaxations and bounds for discrete and nonconvex global optimization. General-purpose LP solvers have traditionally relied on the simplex method~\cite{dantzig2016linear} or interior-point methods~\cite{wright1997primal}. Both frameworks are mature and reliable, but their reliance on sparse matrix factorizations can limit scalability on very large instances. Matrix-free first-order methods avoid these factorizations, relying instead on sparse matrix-vector products and basic vector operations. These operations parallelize efficiently on GPUs and make each iteration inexpensive~\cite{lu2025cupdlp,lu2025overview}. As a result, first-order methods have attracted renewed interest for solving large-scale LPs on modern parallel hardware.

Despite these computational advantages, the difficulty of obtaining high-accuracy solutions within reasonable computation times has historically limited the practical use of first-order methods for LP. \PDLP{} addressed this challenge by combining the primal-dual hybrid gradient method (PDHG) with practical algorithmic enhancements, substantially improving performance on large-scale LPs~\cite{applegate2026pdlp}. Building on this work, \cuPDLP{} adapted the approach to GPUs and demonstrated performance comparable to commercial solvers on standard LP benchmarks~\cite{lu2025cupdlp}. More recently, GPU solvers based on restarted Halpern PDHG iterations, including \cuPDLPx{}, have reported further improvements in computational performance~\cite{lu2025cupdlpx}. Other recent GPU-based first-order LP solvers include \cuPDLPC{}~\cite{lu2023cupdlp}, \HPRLP{}~\cite{chen2026hpr}, and \DPDLP{}~\cite{li2026d}.

A general-purpose LP solver, however, must do more than compute optimal solutions: when an instance is infeasible or unbounded, it should recognize this and return a certificate that can be verified independently. Infeasible LPs arise in branch-and-bound node relaxations, Benders subproblems~\cite{benders1962partitioning}, and model development~\cite{chinneck2008feasibility}. Simplex methods produce Farkas certificates or unbounded rays as a by-product of pivoting~\cite{dantzig2016linear}, and interior-point methods detect primal or dual infeasibility systematically through homogeneous formulations~\cite{ye1994nl,xu1996simplified}. By comparison, the computational studies of recent GPU first-order LP solvers have focused primarily on feasible problems with finite optimal values~\cite{lu2025cupdlp,chen2026hpr,lu2025cupdlpx}.

For PDHG applied directly to the LP saddle-point formulation, the iterates diverge on infeasible or unbounded instances, and certificates of primal or dual infeasibility can be recovered from their asymptotic behavior~\cite{applegate2024infeasibility}; restarted Halpern PDHG accelerates this recovery to a linear rate~\cite{lu2024restarted}. These results provide a theoretical foundation for infeasibility detection with PDHG. In our experience, however, reliable infeasibility and unboundedness detection remains a practical challenge for existing GPU-based PDHG solvers. The limited availability of large-scale infeasible and unbounded LP benchmarks hinders the development and systematic evaluation of such detection schemes.

An alternative is a homogeneous self-dual (HSD) embedding, which encodes primal-dual optimality and infeasibility certificates in a single homogeneous conic system. Suitable nonzero solutions yield either an optimal primal-dual pair or a certificate that can be verified directly. Originally developed for interior-point methods~\cite{ye1994nl,xu1996simplified}, HSD embeddings have since been paired with first-order methods. \SCS{} applies the alternating direction method of multipliers (ADMM) to an HSD embedding of conic programs~\cite{o2016conic}, and \ABIP{} uses ADMM within an HSD interior-point framework for LP~\cite{lin2021admm}. The embedding has also been extended to monotone linear complementarity problems~\cite{o2021operator}.

The HSD system considered here does not directly fit the saddle-point formulation used by standard PDHG. PDHG nevertheless admits a preconditioned proximal-point formulation~\cite{he2012convergence,chambolle2016ergodic,lu2023unified}, with a metric determined by the step sizes and the constraint matrix. The proximal-point framework applies more generally to monotone inclusions~\cite{rockafellar1976monotone}. We use this framework to construct a generalized PDHG operator for the HSD system.

Developing and evaluating a practical solver based on this extension requires addressing the following issues.
(1) Converting two-sided constraint and variable bounds to standard form can introduce additional variables and constraints, enlarging the HSD system.
(2) Representing these bounds directly through box cones avoids explicit bound slacks but couples the coordinates through a shared homogeneous scale. The resulting projections are not coordinatewise separable and require efficient algorithms compatible with a matrix-free implementation.
(3) The homogeneous system always admits the trivial zero solution, so convergence alone does not guarantee an informative outcome.
(4) The \netlib{} infeasible collection is limited in both the number and size of its instances and contains no unbounded LPs. A broader collection of infeasible and unbounded LPs, including larger instances, is therefore needed for a systematic evaluation of certificate recovery.

To address these issues, we make the following contributions.
\begin{enumerate}
    \item We derive a compact HSD formulation for LPs with two-sided constraint and variable bounds. It represents the bounds directly through box cones, without explicit bound slacks, and encodes both primal-dual optimal solutions and infeasibility certificates.

    \item We propose a generalized PDHG operator for this formulation and develop a restarted reflected Halpern method based on it. The two box-cone projections required by the operator reduce to one-dimensional root-finding problems, which we solve using a safeguarded Newton method.

    \item We establish a sufficient initialization condition that keeps the iterates uniformly bounded away from zero. Under this condition, we further prove convergence to a nonzero HSD solution and derive convergence rates for fixed-frequency restarting.

    \item We implement the proposed method in \HPDLP{}, combining GPU computation with practical algorithmic enhancements. Its iterations use only sparse matrix-vector products, vector operations, and reductions. Through numerical comparisons with CPU and GPU solvers, we assess its strengths and limitations on feasible LP benchmarks and on infeasible and unbounded instances.

    \item We introduce \MIPLIBInfeas{}, a benchmark of 229 infeasible and 105 unbounded LPs derived from \MIPLIB{} 2017~\cite{gleixner2021miplib}. The instances range from approximately $10^5$ to $5.9\times 10^7$ nonzeros, supporting the evaluation of certificate recovery at scale.
\end{enumerate}

The primary goal of this paper is to develop a generalized PDHG framework for solving LPs via HSD formulations and to understand and address the methodological and computational challenges of making this approach practical. The numerical comparisons are intended to support this investigation, rather than to establish the superiority of \HPDLP{} over existing LP approaches.

\subsection{Related literature}
\label{subsec:work}

Our work draws on three lines of research: GPU-accelerated first-order methods for LP, infeasibility detection with first-order methods, and algorithms based on homogeneous embeddings.

\paragraph{First-order methods and GPU solvers for LP.}
PDHG is a primal-dual splitting method developed for structured saddle-point problems in imaging~\cite{chambolle2011first,esser2010general}. \PDLP{} combines PDHG with presolve, preconditioning, adaptive step sizes, primal-weight update and adaptive restarts to solve large-scale LPs~\cite{applegate2026pdlp}. 
Restarting~\cite{applegate2023faster} plays a particularly important role in improving both the computational performance and the theoretical convergence guarantees in \PDLP{}. By exploiting the sharpness structure of LPs, suitable restart schemes achieve linear convergence, improving on the sublinear rates established by classical PDHG analyses without strong convexity assumptions~\cite{applegate2023faster}. \cuPDLP{}, implemented in Julia, adapted this framework to GPUs and demonstrated performance comparable to commercial LP solvers on standard benchmarks~\cite{lu2025cupdlp}, which helped stimulate broader interest in GPU-based first-order methods for LP.
A subsequent C implementation, \cuPDLPC{}~\cite{lu2023cupdlp}, made this approach more accessible to practitioners. The later work on restarted Halpern PDHG showed accelerated linear convergence guarantees for both feasible and infeasible LPs~\cite{lu2024restarted}. Building on this method, \cuPDLPx{} provides a modular C implementation with further practical enhancements~\cite{lu2025cupdlpx}. In parallel, \HPRLP{} introduces a semi-proximal Halpern Peaceman--Rachford method and incorporates restarting as a practical heuristic~\cite{chen2026hpr}.
PDLP-based methods with GPU support are also now available in NVIDIA's \cuOpt{}~\cite{blin2024accelerate}, \Gurobi{}~\cite{gurobi2026}, \HiGHS{}~\cite{huangfu2018parallelizing}, \COPT{}~\cite{ge2022copt}, and \Xpress{}~\cite{fico2026xpress}. Beyond LP, GPU-based first-order methods have also been developed for convex quadratic programming~\cite{schubiger2020gpu,lu2026practical,huang2025restarted,chen2025hpr}, second-order cone and other conic programs~\cite{lin2025practical,li2026gpu}, and semidefinite programming~\cite{han2024accelerating,aguirre2025cuhallar,ding2025new}. A broader overview of first-order methods for LP and their GPU implementations is given in~\cite{lu2025overview}.


\paragraph{Infeasibility detection with first-order methods.}
For operator-splitting methods applied directly to primal-dual formulations, infeasibility can be detected from the asymptotic displacement of the iterates. For a class of convex problems with linear or quadratic objectives and conic constraints, successive ADMM iterate differences converge, and their nonzero limits yield infeasibility certificates~\cite{banjac2019infeasibility}. This mechanism is implemented in \OSQP{} for quadratic programs~\cite{stellato2020osqp} and in \COSMO{} for conic programs~\cite{garstka2021cosmo}. Related Douglas--Rachford methods provide diagnostics for infeasibility, unboundedness, and certain pathological cases in conic programming~\cite{liu2017new}. For nonexpansive fixed-point iterations, accelerated sublinear rates for approximating the infimal displacement vector, together with matching lower bounds, were established in~\cite{park2023accelerated}. For PDHG applied directly to the LP saddle-point formulation, Applegate et al.~\cite{applegate2024infeasibility} show that the infimal displacement vector yields primal and dual infeasibility certificates whenever they exist, and establish an improved sublinear rate for the convergence of normalized iterates and normalized averages to this vector. Under a nondegeneracy assumption, successive iterate differences further exhibit eventual linear convergence to this vector~\cite{applegate2024infeasibility}. Subsequent work on restarted Halpern PDHG achieves accelerated eventual linear convergence for certificate recovery without the non-degeneracy assumption~\cite{lu2024restarted}. The displacement-range analysis has also been extended to conic programming, yielding infeasibility results for quadratic programs and a subclass of second-order cone programs~\cite{jiang2023range}. An HSD embedding instead targets a solvable augmented system, allowing a candidate certificate to be recovered from an HSD solution and verified directly.

\paragraph{Homogeneous embeddings and solution methods.}
Homogeneous self-dual embeddings were developed for interior-point LP algorithms~\cite{ye1994nl,xu1996simplified}. The framework was subsequently extended to monotone complementarity problems~\cite{ye1997homogeneous,andersen1999homogeneous} and conic programming~\cite{luo2000conic}, and the recovery of infeasibility certificates from homogeneous models was analyzed in~\cite{andersen2000primal}. Homogeneous embeddings underlie interior-point solvers including \SeDuMi{}~\cite{sturm1999using}, \ECOS{}~\cite{domahidi2013ecos}, \MOSEK{}~\cite{andersen2000mosek}, and \Clarabel{}~\cite{goulart2026clarabel}. \CuClarabel{} provides a GPU implementation of \Clarabel{}~\cite{chen2026cuclarabel}. \Gurobi{} also offers a homogeneous barrier option~\cite{gurobi2026}. \SCS{} applies ADMM to an HSD embedding of conic programs with linear objectives~\cite{o2016conic}, and its extension to quadratic objectives uses Douglas--Rachford splitting on a homogeneous embedding of a monotone linear complementarity problem~\cite{o2021operator}. \ABIP{} and \ABIPp{} use ADMM to approximately solve the barrier subproblems within an interior-point framework based on an HSD formulation~\cite{lin2021admm,deng2022enhanced}. All of these solvers require solving linear systems within their iterations. In contrast, our method applies restarted reflected Halpern iteration to a generalized PDHG operator for the HSD embedding, yielding a GPU-based first-order LP solver that avoids linear system solves.

\section{Preliminaries}
\label{sec:prelim}

In this section, we introduce the notation used throughout the paper (\autoref{subsec:notat}), present the LP formulation and its exact certificates of primal and dual infeasibility (\autoref{subsec:lpcert}), and review PDHG and the Halpern iteration (\autoref{subsec:pdhg}).

\subsection{Notation}
\label{subsec:notat}

Let $\mathbb{R}$ denote the set of real numbers and let $\mathbb{R}_{+}$ and $\mathbb{R}_{-}$ denote the sets of nonnegative and nonpositive real numbers. For a vector $v\in\mathbb{R}^{d}$, its positive and negative parts are denoted by $v^{+}$ and $v^{-}$, where $(v^{+})_{i}:=\max\{v_i,0\}$ and $(v^{-})_{i}:=\max\{-v_i,0\}$. We use $\langle u,v\rangle:=u^{\top}v$ for the Euclidean inner product and $\|\cdot\|_2$ for the Euclidean norm of a vector or the spectral norm of a matrix. For a vector $v\in\mathbb R^d$, we define $\|v\|_\infty:=\max_{1\leq i\leq d}|v_i|$. For bound vectors that may contain infinite entries, we use the convention
\[
    \|(\ell,u)\|_2^2
    :=
    \sum_{i:\,u_i<+\infty}u_i^2
    +
    \sum_{i:\,\ell_i>-\infty,\ \ell_i\neq u_i}\ell_i^2.
\]
Thus, the two finite endpoints of a proper interval are both included, whereas a finite equality bound is counted only once; infinite entries are omitted. For a symmetric positive-definite matrix $M\in\mathbb{R}^{d\times d}$, define
\[
    \langle u,v\rangle_M:=\langle u,Mv\rangle,
    \qquad
    \|u\|_M:=\sqrt{\langle u,u\rangle_M}.
\]
Given a nonempty closed convex set $\mathcal C\subseteq\mathbb{R}^d$, its Euclidean projection and $M$-metric projection are denoted by
\[
    \Pi_{\mathcal C}(x) :=\argmin_{w\in\mathcal C}\|x-w\|_2,
    \qquad
    \Pi_{\mathcal C}^{M}(x) :=\argmin_{w\in\mathcal C}\|x-w\|_M.
\]
The corresponding distance is $\operatorname{dist}_M(x,\mathcal C):=\min_{w\in\mathcal C}\|x-w\|_M$. For an operator $T$, we write $\operatorname{Fix}(T):=\{z: T(z)=z\}$ for its set of fixed points.

\subsection{Linear programming and infeasibility certificates}
\label{subsec:lpcert}

\HPDLP{} solves LPs in the following bound form, with two-sided bounds on both constraints and variables, as used in \PDLP{}~\cite{applegate2026pdlp} and \cuPDLPx{}~\cite{lu2025cupdlpx}:
\begin{equation}
\begin{aligned}
    \min_{x\in\mathbb{R}^n}\quad & c^{\top}x \\
    \text{subject to}\quad
        & \ell_c\leq Ax\leq u_c,\\
        & \ell_v\leq x\leq u_v,
\end{aligned}
\label{eq:primal-lp}
\end{equation}
where $A\in\mathbb{R}^{m\times n}$, $c\in\mathbb{R}^n$, $\ell_c\in(\mathbb{R}\cup\{-\infty\})^m$, $u_c\in(\mathbb{R}\cup\{+\infty\})^m$, $\ell_v\in(\mathbb{R}\cup\{-\infty\})^n$, and $u_v\in(\mathbb{R}\cup\{+\infty\})^n$. Throughout, we assume $\ell_c\leq u_c$ and $\ell_v\leq u_v$.  For convenience, let
\[
    \mathcal{S}:=[\ell_c,u_c],
    \qquad
    \mathcal{X}:=[\ell_v,u_v].
\]

We handle two-sided bounds through the support function of a box: for $\ell\in(\mathbb{R}\cup\{-\infty\})^d$ and $u\in(\mathbb{R}\cup\{+\infty\})^d$ satisfying $\ell\leq u$, define
\begin{equation}\label{eq:p-func}
    p(y;\ell,u)
    :=\sup_{\ell\leq x\leq u}y^\top x
    =u^\top y^{+}-\ell^\top y^{-}.
\end{equation}
Observe that $p(-y;\ell,u)$ is finite if and only if each coordinate $y_{i}$ carries the sign dictated by which of the two bounds are finite. For the constraint bounds $(\ell_c,u_c)$ and the variable bounds $(\ell_v,u_v)$, these sign patterns form the Cartesian products $\mathcal Y\subseteq\mathbb R^m$ and $\mathcal R\subseteq\mathbb R^n$ with coordinate sets
\[
\begin{aligned}
\mathcal Y_i&:=
\begin{cases}
    \{0\}, & \ell_{c,i}=-\infty,\ u_{c,i}=+\infty,\\
    \mathbb R_-, & \ell_{c,i}=-\infty,\ u_{c,i}\in\mathbb R,\\
    \mathbb R_+, & \ell_{c,i}\in\mathbb R,\ u_{c,i}=+\infty,\\
    \mathbb R, & \ell_{c,i},u_{c,i}\in\mathbb R,
\end{cases}
&\quad
\mathcal R_j&:=
\begin{cases}
    \{0\}, & \ell_{v,j}=-\infty,\ u_{v,j}=+\infty,\\
    \mathbb R_-, & \ell_{v,j}=-\infty,\ u_{v,j}\in\mathbb R,\\
    \mathbb R_+, & \ell_{v,j}\in\mathbb R,\ u_{v,j}=+\infty,\\
    \mathbb R, & \ell_{v,j},u_{v,j}\in\mathbb R,
\end{cases}
\end{aligned}
\]
for $i=1,\ldots,m$ and $j=1,\ldots,n$, which are the sign constraints on the dual variables.

Using \eqref{eq:p-func}, \eqref{eq:primal-lp} can be expressed as the convex-concave saddle-point problem
\begin{equation}
    \min_{x\in\mathcal X}\max_{y\in\mathcal Y}
    \quad
    L(x,y)
    :=
    c^\top x
    -y^\top Ax
    -p(-y;\ell_c,u_c).
    \label{eq:lp-saddle}
\end{equation}
The corresponding dual is obtained by minimizing $L(\cdot,y)$ over $\mathcal X$. Introducing $r:=c-A^\top y$ and using $\min_{x\in\mathcal X}r^\top x = - p(-r;\ell_v,u_v)$,
we obtain
\begin{equation}
\begin{aligned}
    \max_{y\in\mathcal Y,\,r\in\mathcal R}\quad&
    -p(-y;\ell_c,u_c)-p(-r;\ell_v,u_v)\\
    \text{subject to}\quad&
    A^\top y+r=c.
\end{aligned}
    \label{eq:dual-lp}
\end{equation}
Thus, at every dual-feasible point, $r=c-A^\top y$ is the reduced-cost vector.

We next recall exact certificates of primal and dual infeasibility. By Farkas' lemma~\cite{boyd2004convex}, \eqref{eq:primal-lp} is primal infeasible if and only if there exists $(y,r)\in\mathcal Y\times\mathcal R$ such that
\begin{equation}
    A^{\top}y+r=0,
    \qquad
    -p(-y;\ell_c,u_c)-p(-r;\ell_v,u_v)>0.
    \label{eq:primal-infeas-cert}
\end{equation}

For a nonempty closed convex set $\mathcal C$, its recession cone is
\[
    \operatorname{rec}(\mathcal C)
    :=\{d:z+td\in\mathcal C\ \text{for every }z\in\mathcal C
    \text{ and }t\geq0\}.
\]
Applying Farkas' lemma to the dual system, \eqref{eq:dual-lp} is dual infeasible if and only if there exists $d\in\mathbb{R}^n$ such that \begin{equation}
    d\in\operatorname{rec}(\mathcal{X}),
    \qquad
    Ad\in\operatorname{rec}(\mathcal{S}),
    \qquad
    c^{\top}d<0.
    \label{eq:dual-infeas-cert}
\end{equation}
If \eqref{eq:primal-lp} is feasible, any such vector $d$ is an improving recession direction and hence certifies that the primal problem is unbounded below.

\subsection{PDHG and the restarted reflected Halpern iteration}
\label{subsec:pdhg}

PDHG~\cite{esser2010general,chambolle2011first} is the base algorithm of \PDLP{}~\cite{applegate2026pdlp}, \cuPDLP{}~\cite{lu2025cupdlp}, \cuPDLPC{}~\cite{lu2023cupdlp}, and \cuPDLPx{}~\cite{lu2025cupdlpx}. Applied to \eqref{eq:lp-saddle}, one PDHG update $z^{k+1}=T_{\mathrm{LP}}(z^k)$, where $z^k:=(x^k,y^k)$, is given by
\begin{subequations}
\label{eq:lp-pdhg}
\begin{align}
    x^{k+1}
    &=\Pi_{\mathcal X}\bigl(x^k-\eta(c-A^{\top}y^k)\bigr),
        \label{eq:lp-pdhg-primal}\\
    y^{k+1}
    &=y^k-\sigma A(2x^{k+1}-x^k)
      -\sigma\Pi_{-\mathcal S}
       \bigl(\sigma^{-1}y^k-A(2x^{k+1}-x^k)\bigr),
        \label{eq:lp-pdhg-dual}
\end{align}
\end{subequations}
where $\eta,\sigma>0$ are the primal and dual step sizes, and $-\mathcal S=[-u_c,-\ell_c]$.

The optimality conditions of \eqref{eq:lp-saddle} can be written as the monotone inclusion $0\in F(z)$ with
\[
    F(x,y):=\bigl(c-A^{\top}y+N_{\mathcal X}(x),\; Ax+\partial p(-y;\ell_c,u_c)\bigr),
\]
where $N_{\mathcal X}$ is the normal cone of $\mathcal X$ and $\partial p(-y;\ell_c,u_c)$ denotes the subdifferential of $y\mapsto p(-y;\ell_c,u_c)$. Following the proximal-point interpretation of PDHG~\cite{he2012convergence,lu2023unified}, the update $z^{k+1}=T_{\mathrm{LP}}(z^k)$ is the unique solution of
\begin{equation}
0\in F(z^{k+1})+P_{\mathrm{LP}}(z^{k+1}-z^k),
\qquad
P_{\mathrm{LP}}
:=
\begin{bmatrix}
\eta^{-1}I_n & A^{\top}\\
A & \sigma^{-1}I_m
\end{bmatrix}.
\label{eq:lp-pdhg-ppa}
\end{equation}
When $\eta\sigma\|A\|_2^2<1$, $P_{\mathrm{LP}}$ is positive definite, and \eqref{eq:lp-pdhg-ppa} identifies PDHG as a preconditioned proximal-point step for $F$ in the $P_{\mathrm{LP}}$-metric~\cite{chambolle2016ergodic}. The fixed points of $T_{\mathrm{LP}}$ are exactly the saddle points of \eqref{eq:lp-saddle}, and monotonicity of $F$ implies that $T_{\mathrm{LP}}$ is firmly nonexpansive in the $P_{\mathrm{LP}}$-norm.

These operator properties provide the basis for accelerating PDHG through Halpern iteration, reflection, and restarting~\cite{lu2024restarted,lu2025cupdlpx}. We briefly review these ingredients below in terms of the operator $T_{\mathrm{LP}}$.

\paragraph{Halpern iteration}
Given an anchor $z^0$, the Halpern iteration~\cite{halpern1967fixed, lu2024restarted} takes the form
\begin{equation}
    z^{k+1}
    =\frac{k+1}{k+2}T_{\mathrm{LP}}(z^k)+\frac{1}{k+2}z^0.
    \label{eq:halpern-basic}
\end{equation}
If $T_{\mathrm{LP}}$ is nonexpansive in the $P_{\mathrm{LP}}$-norm and has a fixed point, the Halpern iterates satisfy $\|z^k-T_{\mathrm{LP}}(z^k)\|_{P_{\mathrm{LP}}}=\mathcal O(1/k)$~\cite{lieder2021convergence}.

\paragraph{Reflection}
Since $T_{\mathrm{LP}}$ is firmly nonexpansive in the $P_{\mathrm{LP}}$-norm, its reflection $2T_{\mathrm{LP}}-I$, where $I$ denotes the identity mapping, is nonexpansive in the same norm and has the same fixed points as $T_{\mathrm{LP}}$~\cite{bauschke2017convex}. Applying the Halpern iteration to $2T_{\mathrm{LP}}-I$ gives the reflected Halpern iteration:
\begin{equation}
    z^{k+1}
    =\frac{k+1}{k+2}\bigl(2T_{\mathrm{LP}}(z^k)-z^k\bigr)
     +\frac{1}{k+2}z^0.
    \label{eq:reflected-halpern}
\end{equation}

\paragraph{Restarting}
The restarted reflected Halpern iteration~\cite{lu2024restarted} resets the anchor when a restart condition is met. Let $n$ and $t$ denote the outer and inner iteration counters. Within outer iteration $n$, it generates
\begin{equation}
    z^{n,t+1}
    =
    \frac{t+1}{t+2}
    \bigl(2T_{\mathrm{LP}}(z^{n,t})-z^{n,t}\bigr)
    +\frac{1}{t+2}z^{n,0}.
    \label{eq:restarted-reflected-halpern}
\end{equation}
In the basic restarted iteration, if a restart is triggered after inner iteration $t_n$, the next outer iteration starts from the corresponding fixed-point step,
\begin{equation}
    z^{n+1,0}:=T_{\mathrm{LP}}(z^{n,t_n}).
    \label{eq:restart-point}
\end{equation}

\HPDLP{} applies this restarted reflected Halpern iteration to the generalized PDHG operator defined in \autoref{sec:algo} for the box-cone homogeneous self-dual formulation derived in \autoref{sec:hsde}.

\section{A Box-Cone Homogeneous Self-Dual Formulation}
\label{sec:hsde}

When an LP is infeasible or unbounded, its Karush-Kuhn-Tucker (KKT) conditions have no solution, and an iterative method applied to them has nothing to converge to. A homogeneous self-dual embedding resolves this by homogenizing the KKT conditions and adding a certificate branch: the resulting conic system admits, for every LP, a solution encoding either primal-dual optimality or an infeasibility certificate.

The form of the embedding matters for a first-order method, since it determines the variables, the cones to project onto, and the matrix-vector products in each iteration. Existing embeddings are formulated for equality-constrained conic programs~\cite{ye1994nl,xu1996simplified,luo2000conic,o2016conic,goulart2026clarabel} or for complementarity problems~\cite{andersen1999homogeneous,o2021operator}, with the cone constraints imposed on the primal variable or on constraint slacks. In particular, the box cone of~\cite{o2021operator} bounds a slack vector relative to a scale coordinate. Here, we impose box cones on $(\tau,x)$ and $(\tau,Ax)$, with the homogeneous variable $\tau$ as the scale coordinate, which represents the two-sided bounds in \eqref{eq:primal-lp} without slack variables.

The derivation proceeds in four steps. We first homogenize the KKT conditions so that certificates appear as solutions with zero scale (\autoref{subsec:kkt}), represent the scaled two-sided bounds as box cones (\autoref{subsec:boxcone}), eliminate the reduced costs and the gap variable to obtain a skew-symmetric conic system (\autoref{subsec:conic}), and characterize its solutions (\autoref{subsec:solu}).

\subsection{Homogenizing the KKT conditions}
\label{subsec:kkt}

At a primal-feasible point $x$ and a dual-feasible point $(y,r)$, the difference between the primal and dual objectives is
\[
    c^{\top}x+p(-y;\ell_c,u_c)+p(-r;\ell_v,u_v).
\]
It is nonnegative by weak duality and vanishes at optimality, so the KKT conditions for the primal-dual pair \eqref{eq:primal-lp} and \eqref{eq:dual-lp} can be written as
\begin{equation}
\begin{aligned}
    &x\in\mathcal X,\qquad Ax\in\mathcal S,\\
    &y\in\mathcal Y,\qquad r\in\mathcal R,\qquad
      A^{\top}y+r=c,\\
    &c^{\top}x+p(-y;\ell_c,u_c)+p(-r;\ell_v,u_v)=0.
\end{aligned}
\label{eq:kkt-box}
\end{equation}

Since \eqref{eq:kkt-box} has no solution when the LP is infeasible or unbounded, we embed it in a system that also admits the certificates \eqref{eq:primal-infeas-cert} and \eqref{eq:dual-infeas-cert} as solutions. To this end, we introduce a scale $\tau>0$ and require the normalized point $(x/\tau,y/\tau,r/\tau)$ to satisfy \eqref{eq:kkt-box}. Since $\mathcal Y$ and $\mathcal R$ are cones and the support functions are positively homogeneous, this is equivalent to
\begin{equation}
\begin{aligned}
    &x\in\tau\mathcal X,\qquad Ax\in\tau\mathcal S,\\
    &y\in\mathcal Y,\qquad r\in\mathcal R,\qquad
      A^{\top}y+r=c\tau,\\
    &c^{\top}x+p(-y;\ell_c,u_c)+p(-r;\ell_v,u_v)=0.
\end{aligned}
\label{eq:scaled-kkt-box}
\end{equation}
Setting $\tau=0$ in \eqref{eq:scaled-kkt-box}, the dual conditions give the feasibility and stationarity conditions for the primal infeasibility certificate \eqref{eq:primal-infeas-cert}; the primal conditions give the feasibility and stationarity conditions for \eqref{eq:dual-infeas-cert} once the scaled boxes are extended to closed cones in \autoref{subsec:boxcone}. A certificate of either kind, however, violates the gap equation: it makes $c^{\top}x+p(-y;\ell_c,u_c)+p(-r;\ell_v,u_v)$ strictly negative rather than zero. We therefore relax the gap equation by a nonnegative scalar,
\begin{equation}
    c^{\top}x+p(-y;\ell_c,u_c)+p(-r;\ell_v,u_v)+\kappa=0,
    \qquad \kappa\geq0.
\label{eq:hsd-gap}
\end{equation}
For $\tau>0$, weak duality forces $\kappa=0$, so the augmentation leaves the KKT conditions unchanged. For $\tau=0$, every solution with $\kappa>0$ yields a certificate, and every certificate yields such a solution (\autoref{subsec:solu}).

\subsection{Box-cone representation}
\label{subsec:boxcone}

At $\tau=0$, the scaled boxes $\tau\mathcal X$ and $\tau\mathcal S$ in \eqref{eq:scaled-kkt-box} reduce to $\{0\}$, which excludes the recession directions required by \eqref{eq:dual-infeas-cert}. We therefore replace them by closed cones in $(\tau,x)$. For vectors $\ell\in(\mathbb{R}\cup\{-\infty\})^d$ and $u\in(\mathbb{R}\cup\{+\infty\})^d$ satisfying $\ell\leq u$, define the \emph{box cone}
\begin{equation}
    \mathcal K[\ell,u]
    :=\operatorname{cl}\bigl\{(\tau,x)\in\mathbb{R}\times\mathbb{R}^d:
       \tau>0,\ \ell\tau\leq x\leq u\tau\bigr\},
\label{eq:box-cone}
\end{equation}
where $\operatorname{cl}$ denotes set closure. For $\tau>0$, $(\tau,x)\in\mathcal K[\ell,u]$ means $x/\tau\in[\ell,u]$, and at $\tau=0$,
\begin{equation}
    (0,x)\in\mathcal K[\ell,u]
    \quad\Longleftrightarrow\quad
    x\in\operatorname{rec}([\ell,u])
\label{eq:box-cone-recession}
\end{equation}
by~\cite{rockafellar1970convex}. For finite bounds, $\mathcal K[\ell,u]=\{(\tau,x):\tau\geq0,\ \ell\tau\leq x\leq u\tau\}$ is the box cone of~\cite{o2021operator}, and taking the closure adds only the origin; when some bounds are infinite, the closure may also add nonzero recession directions.

Applying this construction to the variable and constraint boxes, set $\mathcal K_v:=\mathcal K[\ell_v,u_v]$ and $\mathcal K_c:=\mathcal K[\ell_c,u_c]$. The primal conditions in \eqref{eq:scaled-kkt-box} then extend to
\begin{equation}
    (\tau,x)\in\mathcal K_v,
    \qquad
    (\tau,Ax)\in\mathcal K_c.
\label{eq:primal-box-cones}
\end{equation}
The remaining conditions, dual feasibility and the gap condition, will be expressed in \autoref{subsec:conic} through the dual cones $\mathcal K_v^*$ and $\mathcal K_c^*$, which the support function \eqref{eq:p-func} characterizes as follows.

\begin{lem}
\label{lem:dual-box-cone}
For any $\ell$ and $u$ as above,
\begin{equation}
    \mathcal K[\ell,u]^*
    =\bigl\{(\alpha,y)\in\mathbb{R}\times\mathbb{R}^d:
      \alpha\geq p(-y;\ell,u)\bigr\}.
\label{eq:dual-box-cone}
\end{equation}
\end{lem}

\begin{proof}
By definition, $(\alpha,y)\in\mathcal K[\ell,u]^*$ if and only if
\[
    \alpha\tau+y^\top x\geq0
    \qquad
    \text{for every }(\tau,x)\in\mathcal K[\ell,u].
\]
By continuity, it suffices to check the inequality before taking the closure in \eqref{eq:box-cone}. Writing $x=\tau w$ with $\tau>0$ and $w\in[\ell,u]$, the condition becomes
\[
    \tau\bigl(\alpha+y^\top w\bigr)\geq0
    \qquad
    \text{for every }\tau>0\text{ and }w\in[\ell,u].
\]
Since $\tau>0$, this holds if and only if
\[
    \alpha\geq(-y)^\top w
    \qquad
    \text{for every }w\in[\ell,u],
\]
or equivalently,
\[
    \alpha\geq\sup_{w\in[\ell,u]}(-y)^\top w
    =p(-y;\ell,u).
\]
\qed
\end{proof}

\subsection{Skew-symmetric conic form}
\label{subsec:conic}

The primal constraints are now represented by \eqref{eq:primal-box-cones}. It remains to express dual feasibility and the gap condition \eqref{eq:hsd-gap} in conic form. Rather than carrying $r$ and $\kappa$ as variables, we treat them as functions of $(\tau,x,y)$:
\begin{equation}
\begin{aligned}
    r&:=c\tau-A^{\top}y,\\
    \kappa&:=-c^{\top}x
      -p(-y;\ell_c,u_c)-p(-r;\ell_v,u_v).
\end{aligned}
\label{eq:hsd-recovery}
\end{equation}
With $r$ defined by \eqref{eq:hsd-recovery}, the equation $A^{\top}y+r=c\tau$ holds identically, and the support-function values are finite exactly when $y\in\mathcal Y$ and $r\in\mathcal R$. Hence, dual feasibility together with $\kappa\geq0$ is equivalent to
\begin{equation}
    p(-y;\ell_c,u_c)+p(-r;\ell_v,u_v)
    \leq-c^{\top}x,
\label{eq:hsd-gap-inequality}
\end{equation}
where \eqref{eq:hsd-gap-inequality} is read as violated whenever its left-hand side is $+\infty$.
Introducing a scalar $\alpha$ that bounds the first support function, \eqref{eq:hsd-gap-inequality} holds if and only if there exists $\alpha\in\mathbb R$ with
\begin{equation}
\begin{aligned}
    \alpha
    &\geq p(-y;\ell_c,u_c),\\
    -c^{\top}x-\alpha
    &\geq p(-r;\ell_v,u_v),
\end{aligned}
\label{eq:alpha-split}
\end{equation}
since adding the two inequalities gives \eqref{eq:hsd-gap-inequality}, and conversely $\alpha=p(-y;\ell_c,u_c)$ satisfies both. By \eqref{eq:dual-box-cone}, the two inequalities in \eqref{eq:alpha-split} are equivalent to
\[
    (\alpha,y)\in\mathcal K_c^*,
    \qquad
    (-c^{\top}x-\alpha,r)\in\mathcal K_v^*.
\]
Combining these constraints with \eqref{eq:primal-box-cones} gives
\begin{equation}
\begin{aligned}
    &(\tau,x)\in\mathcal K_v,
    &&(\tau,Ax)\in\mathcal K_c,\\
    &(\alpha,y)\in\mathcal K_c^*,
    &&(-c^{\top}x-\alpha,c\tau-A^{\top}y)\in\mathcal K_v^*.
\end{aligned}
\label{eq:four-cone-system}
\end{equation}

To write the four cone constraints \eqref{eq:four-cone-system} compactly, let
$z:=(\tau,x,\alpha,y)$ and define
\begin{equation}
    Q:=
    \begin{bmatrix}
        0 & -c^{\top} & -1 & 0\\
        c & 0 & 0 & -A^{\top}\\
        1 & 0 & 0 & 0\\
        0 & A & 0 & 0
    \end{bmatrix},
    \qquad
    \mathcal K:=\mathcal K_v\times\mathcal K_c^*.
\label{eq:hsd-QK}
\end{equation}
Then $\mathcal K^*=\mathcal K_v^*\times\mathcal K_c$ and
\[
    Qz=
    \begin{bmatrix}
        -c^{\top}x-\alpha\\
        c\tau-A^{\top}y\\
        \tau\\
        Ax
    \end{bmatrix}.
\]
Therefore, \eqref{eq:four-cone-system} is equivalent to the compact conic system
\begin{equation}
    z\in\mathcal K,
    \qquad
    Qz\in\mathcal K^*.
\label{eq:box-hsd}
\end{equation}

Since $Q$ is skew-symmetric, \eqref{eq:box-hsd} has the standard form of a homogeneous self-dual embedding~\cite{o2016conic}. Its solution set, denoted $\mathcal Z^\star$, is the intersection of $\mathcal K$ with the closed convex cone $\{z:Qz\in\mathcal K^*\}$, hence a closed convex cone containing the origin.

As intended, $z$ consists only of $(\tau,x,\alpha,y)$, $\mathcal K$ is a product of a box cone and a dual box cone, and $Q$ involves the problem data only through $c$, $A$, and $A^{\top}$. Therefore, the PDHG iteration of \autoref{sec:algo} requires the same matrix-vector products as on the original LP, plus two box-cone projections.

\subsection{Solution characterization}
\label{subsec:solu}
As in standard homogeneous models~\cite{andersen2000primal,o2016conic}, solutions with $\tau>0$ yield primal-dual optimal solutions, whereas solutions with $\tau=0$ and $\kappa>0$ yield infeasibility certificates. The following proposition makes this precise for \eqref{eq:box-hsd}.

\begin{prop}
\label{prop:hsd-solution-characterization}
Let $z=(\tau,x,\alpha,y)$ satisfy \eqref{eq:box-hsd}, and let $r$ and $\kappa$ be given by \eqref{eq:hsd-recovery}. Then
\[
    \tau\geq0,
    \qquad
    \kappa\geq0,
    \qquad
    \tau\kappa=0.
\]
If $\tau>0$, then $(x/\tau,y/\tau,r/\tau)$ is a primal-dual optimal solution of \eqref{eq:primal-lp} and \eqref{eq:dual-lp}. If $\kappa>0$, then at least one of
\[
    -p(-y;\ell_c,u_c)-p(-r;\ell_v,u_v)>0,
    \qquad
    c^{\top}x<0
\]
holds. In the first case, $(y,r)$ is a certificate of primal infeasibility in the sense of \eqref{eq:primal-infeas-cert}; in the second, $x$ is a certificate of dual infeasibility in the sense of \eqref{eq:dual-infeas-cert}.
\end{prop}

\begin{proof}
Since $(\tau,x)\in\mathcal K_v$, we have $\tau\geq0$. By \eqref{eq:dual-box-cone}, the two dual-cone constraints in \eqref{eq:four-cone-system} are exactly the inequalities \eqref{eq:alpha-split}. In particular, both support-function values in \eqref{eq:hsd-recovery} are finite, and adding the inequalities gives $\kappa\geq0$.

For $\tau\kappa=0$, note that $(p(-r;\ell_v,u_v),r)\in\mathcal K_v^*$ and $(p(-y;\ell_c,u_c),y)\in\mathcal K_c^*$ by \eqref{eq:dual-box-cone}, so pairing them with $(\tau,x)\in\mathcal K_v$ and $(\tau,Ax)\in\mathcal K_c$ gives $x^{\top}r\geq-\tau p(-r;\ell_v,u_v)$ and $y^{\top}Ax\geq-\tau p(-y;\ell_c,u_c)$. Since $Q$ is skew-symmetric,
\begin{align*}
    0
    =\langle z,Qz\rangle
    &=-\tau c^{\top}x+x^{\top}r+y^{\top}Ax \\
    &\geq\tau\bigl[
        -c^{\top}x
        -p(-y;\ell_c,u_c)
        -p(-r;\ell_v,u_v)
    \bigr]
    =\tau\kappa,
\end{align*}
where the terms in $\alpha$ cancel in the first equality. Since $\tau,\kappa\geq0$, we conclude that $\tau\kappa=0$.

If $\tau>0$, dividing the cone constraints in \eqref{eq:four-cone-system} by $\tau$ shows that $(x/\tau,y/\tau,r/\tau)$ is primal-dual feasible. Since $\kappa=0$, positive homogeneity of the support functions gives zero duality gap, so it satisfies the KKT conditions \eqref{eq:kkt-box}.

If $\kappa>0$, then $\tau=0$. The recession identity \eqref{eq:box-cone-recession} gives $x\in\operatorname{rec}(\mathcal X)$ and $Ax\in\operatorname{rec}(\mathcal S)$, while $r=c\tau-A^{\top}y$ gives $A^{\top}y+r=0$. The finiteness established above gives $(y,r)\in\mathcal Y\times\mathcal R$. Finally, \eqref{eq:hsd-recovery} implies that at least one of the two strict inequalities in the proposition holds, yielding the corresponding certificate in \eqref{eq:primal-infeas-cert} or \eqref{eq:dual-infeas-cert}.
\qed
\end{proof}

A solution with $\tau=\kappa=0$ is not conclusive unless it happens to satisfy one of the two certificate conditions. The next result shows that a conclusive solution always exists.

\begin{cor}
\label{cor:informative-hsd-solution}
The homogeneous self-dual formulation \eqref{eq:box-hsd} admits a solution
satisfying $\tau+\kappa>0$.
\end{cor}

\begin{proof}
By LP duality, exactly one of the following holds. If a primal-dual optimal solution $(x,y,r)$ exists, with $r=c-A^{\top}y$, then setting $\tau=1$ and $\alpha=p(-y;\ell_c,u_c)$ gives a solution of \eqref{eq:box-hsd} with $\kappa=0$. If the primal problem is infeasible, a certificate satisfying \eqref{eq:primal-infeas-cert}, together with $\tau=0$, $x=0$, and $\alpha=p(-y;\ell_c,u_c)$, gives a solution with $\kappa>0$. Finally, if the primal problem is feasible and unbounded below, an improving recession direction $d$ satisfying \eqref{eq:dual-infeas-cert}, together with $\tau=0$, $x=d$, $\alpha=0$, and $y=0$, gives a solution with $\kappa=-c^{\top}d>0$.
\qed
\end{proof}

Although \autoref{cor:informative-hsd-solution} guarantees a solution with $\tau+\kappa>0$, the origin also satisfies \eqref{eq:box-hsd}, and an iterative method may converge to it. In \autoref{subsec:op}, we express \eqref{eq:box-hsd} as a variational inequality and define its generalized PDHG operator. \autoref{subsec:init} shows how the initialization keeps the iterates away from this trivial solution.

\section{Halpern PDHG for the Box-Cone HSD Formulation}
\label{sec:algo}

We now develop a restarted reflected Halpern PDHG method for the HSD system \eqref{eq:box-hsd}. We construct the generalized PDHG operator and establish its fixed-point and nonexpansiveness properties (\autoref{subsec:op}), then reduce its box-cone updates to one-dimensional root-finding problems (\autoref{subsec:proj}). For fixed step sizes and the basic restart rule \eqref{eq:restart-point}, we then give an initialization that keeps the iterates uniformly bounded away from the trivial zero solution (\autoref{subsec:init}). Finally, we establish convergence rates for the exact iteration (\autoref{subsec:convergence}).

\subsection{The generalized PDHG operator}
\label{subsec:op}

Since $Q$ is skew-symmetric and $\mathcal K$ is a cone, the HSD system \eqref{eq:box-hsd} is equivalent to the variational inequality
\begin{equation}
    z\in\mathcal K,
    \qquad
    \langle Qz,v-z\rangle\geq0
    \qquad\text{for every }v\in\mathcal K.
    \label{eq:hsd-vi}
\end{equation}
Indeed, skew-symmetry gives $\langle Qz,v-z\rangle=\langle Qz,v\rangle$, so \eqref{eq:hsd-vi} is equivalent to $z\in\mathcal K$ and $Qz\in\mathcal K^*$.

For $u\in\mathcal K$, let
\[
    N_{\mathcal K}(u):=\left\{w:\ \langle w,v-u\rangle\leq0\ \text{for every }v\in\mathcal K\right\}
\]
denote the normal cone of $\mathcal K$ at $u$, and set $N_{\mathcal K}(u):=\emptyset$ for $u\notin\mathcal K$. Then \eqref{eq:hsd-vi} is equivalent to the monotone inclusion
\begin{equation}
    0\in Qz+N_{\mathcal K}(z).
    \label{eq:hsd-monotone-inclusion}
\end{equation}
Since $Q$ is skew-symmetric and $N_{\mathcal K}$ is a normal-cone operator, both are monotone, and so is their sum in \eqref{eq:hsd-monotone-inclusion}.

Under the splitting $(\tau,x)$ and $(\alpha,y)$, \eqref{eq:hsd-monotone-inclusion} is not the optimality condition of a saddle-point problem: the terms $-c^{\top}x$ and $c\tau$ form a skew-symmetric coupling within the primal block. Motivated by the proximal-point interpretation of PDHG~\cite{he2012convergence,lu2023unified,lu2025overview}, we define the iteration for \eqref{eq:hsd-monotone-inclusion} through the proximal-point form \eqref{eq:lp-pdhg-ppa}, with the preconditioner chosen to retain the block structure of PDHG.

Let $\mu,\eta,\sigma>0$ be the step sizes for $\tau$, $x$, and $(\alpha,y)$, and define
\begin{equation}
    P:=
    \begin{bmatrix}
        \mu^{-1} & 0 & -1 & 0\\
        0 & \eta^{-1}I_n & 0 & -A^{\top}\\
        -1 & 0 & \sigma^{-1} & 0\\
        0 & -A & 0 & \sigma^{-1}I_m
    \end{bmatrix}.
    \label{eq:hsd-pdhg-metric}
\end{equation}
The off-diagonal blocks of $P$ are chosen so that the lower-left block of $Q+P$ vanishes, making $Q+P$ block upper triangular under the splitting $(\tau,x)$ and $(\alpha,y)$. We choose the step sizes so that $P\succ0$, or equivalently,
\begin{equation}
    \mu\sigma<1
    \qquad\text{and}\qquad
    \eta\sigma\|A\|_2^2<1.
    \label{eq:hsd-step-condition}
\end{equation}

For any $z\in\mathbb R^{n+m+2}$, define $T(z)$ through the preconditioned proximal-point inclusion
\begin{equation}
    0\in QT(z)+N_{\mathcal K}\bigl(T(z)\bigr)+P\bigl(T(z)-z\bigr).
    \label{eq:hsd-pdhg-inclusion}
\end{equation}
Equivalently, $T(z)\in\mathcal K$ satisfies
\begin{equation}
    \bigl\langle QT(z)+P\bigl(T(z)-z\bigr), v-T(z)\bigr\rangle\geq0
    \qquad\text{for every }v\in\mathcal K.
    \label{eq:hsd-pdhg-proximal}
\end{equation}

For fixed $z$, the mapping $u\mapsto Qu+P(u-z)$ is continuous, and skew-symmetry of $Q$ gives
\[
    \bigl\langle (Q+P)(u-w),u-w\bigr\rangle
    =\|u-w\|_P^2
    \geq \lambda_{\min}(P)\|u-w\|_2^2.
\]
Thus, this mapping is strongly monotone. Since $\mathcal K$ is nonempty, closed, and convex, \eqref{eq:hsd-pdhg-proximal} has a unique solution for every input $z$~\cite{facchinei2003finite}, so $T$ is well defined.

Writing \eqref{eq:hsd-pdhg-proximal} at two inputs $z$ and $w$, taking $v=T(w)$ and $v=T(z)$ respectively, and adding the two inequalities gives
\[
    \|T(z)-T(w)\|_P^2\leq\langle T(z)-T(w),z-w\rangle_P.
\]
Hence, $T$ is firmly nonexpansive in the $P$-norm~\cite{bauschke2017convex}. Moreover, if $T(z)=z$, then \eqref{eq:hsd-pdhg-inclusion} reduces to \eqref{eq:hsd-monotone-inclusion}. Conversely, if $z$ satisfies \eqref{eq:hsd-monotone-inclusion}, then it solves the proximal subproblem with input $z$, and uniqueness gives $T(z)=z$. Thus,
\[
    \operatorname{Fix}(T)=\mathcal Z^\star,
\]
where $\mathcal Z^\star$ denotes the solution set of \eqref{eq:box-hsd}.

Accordingly, \HPDLP{} runs \eqref{eq:restarted-reflected-halpern}--\eqref{eq:restart-point} with $T$. Termination tests are evaluated at the PDHG point $T(z^{n,t})$, and restarts are governed by the fixed-point residual between $z^{n,t}$ and $T(z^{n,t})$ (\autoref{subsec:restrt}).

\paragraph{Computing the PDHG step}
To describe a single evaluation of $T$, write $z^k=(\tau^k,x^k,\alpha^k,y^k)$ and $T(z^k)=(\tau^{k+1},x^{k+1},\alpha^{k+1},y^{k+1})$. Since $Q+P$ is block upper triangular and $\mathcal K=\mathcal K_v\times\mathcal K_c^*$, \eqref{eq:hsd-pdhg-inclusion} can be solved first for $(\alpha,y)$ and then for $(\tau,x)$. With $H:=\operatorname{diag}(\mu^{-1},\eta^{-1}I_n)$, the two updates are
\begin{subequations}
\label{eq:hsd-pdhg}
\begin{align}
    \begin{bmatrix}\alpha^{k+1}\\ y^{k+1}\end{bmatrix}
    &=\Pi_{\mathcal K_c^*}\left(
        \begin{bmatrix}
            \alpha^k-\sigma\tau^k\\
            y^k-\sigma Ax^k
        \end{bmatrix}\right),
    \label{eq:hsd-pdhg-dual}\\
    \begin{bmatrix}\tau^{k+1}\\ x^{k+1}\end{bmatrix}
    &=\Pi_{\mathcal K_v}^{H}\left(
        \begin{bmatrix}
            \tau^k+\mu\bigl(c^{\top}x^{k+1}+2\alpha^{k+1}-\alpha^k\bigr)\\
            x^k-\eta\bigl(c\tau^{k+1}-A^{\top}(2y^{k+1}-y^k)\bigr)
        \end{bmatrix}\right).
    \label{eq:hsd-pdhg-primal}
\end{align}
\end{subequations}
The primal update is implicit: $\tau^{k+1}$ and $x^{k+1}$ are coupled through $c$ and must be determined jointly. Both updates reduce to one-dimensional subproblems, as shown in \autoref{subsec:proj}.

\subsection{Box-cone projections}
\label{subsec:proj}

We first derive the scalar reduction for Euclidean projection onto a box cone, for which a scalar Newton method was noted in~\cite{o2021operator} for finite bounds. We then describe a common scalar solver and its use in the dual and implicit primal updates.

\paragraph{Euclidean projection}
The Euclidean projection of $(\widetilde\tau,\widetilde w)$ onto $\mathcal K[\ell,u]$ is the unique solution of
\begin{equation}
    \min_{(\tau,w)\in\mathcal K[\ell,u]}
    \quad
    \frac12(\tau-\widetilde\tau)^2
    +\frac12\|w-\widetilde w\|_2^2.
    \label{eq:box-cone-projection}
\end{equation}
For fixed $\tau>0$, minimizing over $w$ gives
\begin{equation}
    w_i(\tau)
    =\min\{\max\{\widetilde w_i,\ell_i\tau\},u_i\tau\}.
    \label{eq:box-cone-clipping}
\end{equation}
At $\tau=0$, \eqref{eq:box-cone-recession} gives $w(0)=\Pi_{\operatorname{rec}([\ell,u])}(\widetilde w)$. Eliminating $w$ therefore reduces \eqref{eq:box-cone-projection} to
\begin{equation}
    \min_{\tau\geq0}
    \quad
    \phi(\tau)
    :=\frac12(\tau-\widetilde\tau)^2
      +\frac12\|w(\tau)-\widetilde w\|_2^2.
    \label{eq:box-cone-scalar-projection}
\end{equation}
Let $\mathcal I_{\ell}(\tau)$ and $\mathcal I_u(\tau)$ collect the coordinates that require clipping to the lower and upper bounds, respectively:
\[
    \mathcal I_{\ell}(\tau)
    :=\{i:\ell_i>-\infty,\ \widetilde w_i<\ell_i\tau\},
    \qquad
    \mathcal I_u(\tau)
    :=\{i:u_i<+\infty,\ \widetilde w_i>u_i\tau\}.
\]
Differentiating \eqref{eq:box-cone-scalar-projection} gives
\begin{equation}
    \phi'(\tau)
    =\tau-\widetilde\tau
     +\sum_{i\in\mathcal I_{\ell}(\tau)}
        \ell_i(\ell_i\tau-\widetilde w_i)
     +\sum_{i\in\mathcal I_u(\tau)}
        u_i(u_i\tau-\widetilde w_i),
    \label{eq:box-cone-scalar-derivative}
\end{equation}
where $\phi'(0)$ denotes the right derivative. On intervals with a fixed clipping pattern, $\phi'$ is affine with slope
\begin{equation}
    d(\tau)
    :=1+\sum_{i\in\mathcal I_{\ell}(\tau)}\ell_i^2
        +\sum_{i\in\mathcal I_u(\tau)}u_i^2
    \geq1.
    \label{eq:box-cone-scalar-curvature}
\end{equation}
The affine pieces agree at the breakpoints, so $\phi'$ is continuous. Together with the slope bound $d(\tau)\geq1$, this implies that $\phi'$ is strictly increasing and $\phi$ is $1$-strongly convex. Both $\phi'$ and $d$ can be evaluated using coordinatewise operations and reductions.

\paragraph{Safeguarded Newton method}
We solve the scalar optimality conditions using a safeguarded Newton method, which is also used for the implicit primal update below. To describe the common routine, let $g$ be a continuous piecewise-affine function with slopes bounded below by $d_{\min}>0$. We seek $\tau\geq0$ satisfying $g(\tau)\geq0$ and $\tau g(\tau)=0$. For Euclidean projection, $g=\phi'$, $d_{\min}=1$, and the Newton denominator $d$ is given by \eqref{eq:box-cone-scalar-curvature}, including at clipping breakpoints. If $g(0)\geq0$, the solution is $\tau=0$. Otherwise, the unique positive root lies in $[0,-g(0)/d_{\min}]$, since $g(\tau)\geq g(0)+d_{\min}\tau$.

\autoref{alg:box-cone-newton} maintains this bracket and accepts a Newton step only if it stays within the bracket and its length is at most half that of the previous step. Otherwise, it takes a bisection step. This safeguard guarantees convergence for the piecewise-affine functions considered here. The stopping tests guarantee a relative error of at most $\epsilon_{\mathrm N}$ in the positive root. When available, the solution of the corresponding subproblem in the previous PDHG step provides a warm start. In our implementation, we set $\epsilon_{\mathrm N}=10^{-6}$ and limit each subproblem solve to 100 Newton iterations.

\begin{algorithm}[t]
\caption{Safeguarded Newton method for the scalar subproblems}
\label{alg:box-cone-newton}
\begin{algorithmic}[1]
\REQUIRE Scalar function $g$; Newton denominator $d$; slope lower bound $d_{\min}>0$; relative tolerance $\epsilon_{\mathrm N}>0$; optional warm start $\tau_{\mathrm{ws}}$.
\IF{$g(0)\geq 0$}
  \RETURN $0$
\ENDIF
\STATE $L\gets 0$, $U\gets -g(0)/d_{\min}$, $\Delta\gets +\infty$.
\IF{$\tau_{\mathrm{ws}}$ is available}
  \STATE $\tau\gets \Pi_{[L,U]}(\tau_{\mathrm{ws}})$.
\ELSE
  \STATE $\tau\gets (L+U)/2$.
\ENDIF
\WHILE{true}
  \IF{$|g(\tau)|\leq d_{\min}\frac{\epsilon_{\mathrm N}}{1+\epsilon_{\mathrm N}}\tau$}
    \RETURN $\tau$
  \ENDIF
  \IF{$g(\tau)<0$}
    \STATE $L\gets \tau$.
  \ELSE
    \STATE $U\gets \tau$.
  \ENDIF
  \IF{$U-L\leq 2\epsilon_{\mathrm N}L$}
    \RETURN $(L+U)/2$
  \ENDIF
  \STATE $\tau_{\mathrm{new}}\gets \tau-g(\tau)/d(\tau)$.
  \IF{$\tau_{\mathrm{new}}\notin [L,U]$ \OR $|\tau_{\mathrm{new}}-\tau|>\Delta/2$}
    \STATE $\tau_{\mathrm{new}}\gets (L+U)/2$.
  \ENDIF
  \STATE $\Delta\gets |\tau_{\mathrm{new}}-\tau|$, $\tau\gets \tau_{\mathrm{new}}$.
\ENDWHILE
\end{algorithmic}
\end{algorithm}
For Euclidean projection, the vector component is recovered from \eqref{eq:box-cone-clipping}, with $w(0)$ as defined above.

\paragraph{Dual update}
Let
\[
    \widetilde\alpha:=\alpha^k-\sigma\tau^k,
    \qquad
    \widetilde y:=y^k-\sigma Ax^k.
\]
The dual step \eqref{eq:hsd-pdhg-dual} is $\Pi_{\mathcal K_c^*}(\widetilde\alpha,\widetilde y)$. By the Moreau decomposition~\cite{moreau1962decomposition},
\begin{equation}
    \Pi_{\mathcal K_c^*}(\widetilde\alpha,\widetilde y)
    =(\widetilde\alpha,\widetilde y)
      +\Pi_{\mathcal K_c}(-\widetilde\alpha,-\widetilde y).
    \label{eq:box-cone-moreau}
\end{equation}
Thus, the dual step requires one Euclidean projection onto $\mathcal K_c$, computed by the procedure above.

\paragraph{Primal update}
Let
\begin{equation}
    \widetilde\tau:=\tau^k+\mu(2\alpha^{k+1}-\alpha^k),
    \qquad
    \widetilde x:=x^k+\eta A^{\top}(2y^{k+1}-y^k).
    \label{eq:primal-box-shift}
\end{equation}
The primal step \eqref{eq:hsd-pdhg-primal} becomes
\begin{equation}
    \begin{bmatrix}\tau^{k+1}\\x^{k+1}\end{bmatrix}
    =\Pi_{\mathcal K_v}^{H}\left(
        \begin{bmatrix}
            \widetilde\tau+\mu c^{\top}x^{k+1}\\
            \widetilde x-\eta c\tau^{k+1}
        \end{bmatrix}\right).
    \label{eq:primal-box-rewritten}
\end{equation}
For fixed $\tau>0$, its vector component is
\begin{equation}
    x_i(\tau)
    =\min\{\max\{\widetilde x_i-\eta c_i\tau,\ell_{v,i}\tau\},u_{v,i}\tau\},
    \label{eq:primal-box-clipping}
\end{equation}
with $x(0)=\Pi_{\operatorname{rec}(\mathcal X)}(\widetilde x)$. For this update, $\mathcal I_{\ell}(\tau)$ and $\mathcal I_u(\tau)$ are defined as before with $\widetilde x_i-\eta c_i\tau$ in place of $\widetilde w_i$ and $(\ell_v,u_v)$ in place of $(\ell,u)$.

Eliminating the vector component from the projection optimality conditions gives $\tau\geq0$, $\psi(\tau)\geq0$, and $\tau\psi(\tau)=0$, where
\begin{equation}
\begin{aligned}
    \psi(\tau)
    :={}&\frac{\tau-\widetilde\tau}{\mu}-c^{\top}x(\tau)\\
    &+\sum_{i\in\mathcal I_{\ell}(\tau)}
        \ell_{v,i}\left(\frac{x_i(\tau)-\widetilde x_i}{\eta}+c_i\tau\right)\\
    &+\sum_{i\in\mathcal I_u(\tau)}
        u_{v,i}\left(\frac{x_i(\tau)-\widetilde x_i}{\eta}+c_i\tau\right).
\end{aligned}
    \label{eq:primal-box-scalar-residual}
\end{equation}

As in the Euclidean case, $\psi$ is continuous and piecewise affine. We use the Newton denominator
\[
    d(\tau)
    :=\mu^{-1}
      +\eta\sum_{i\notin\mathcal I_{\ell}(\tau)\cup\mathcal I_u(\tau)}c_i^2
      +\eta^{-1}\sum_{i\in\mathcal I_{\ell}(\tau)}\ell_{v,i}^2
      +\eta^{-1}\sum_{i\in\mathcal I_u(\tau)}u_{v,i}^2.
\]
Away from clipping breakpoints, $d(\tau)$ is the slope of $\psi$, and $d(\tau)\geq\mu^{-1}$. We apply \autoref{alg:box-cone-newton} with $g=\psi$ and $d_{\min}=\mu^{-1}$, giving $\tau^{k+1}=0$ if $\psi(0)\geq0$ and otherwise bracketing the root
in $[0,-\mu\psi(0)]$. Finally, $x^{k+1}=x(\tau^{k+1})$ follows from \eqref{eq:primal-box-clipping}, using $x(0)$ when $\tau^{k+1}=0$.

\subsection{Initialization and avoidance of the trivial solution}
\label{subsec:init}

The HSD system \eqref{eq:box-hsd} always admits $z=0$. The following proposition gives an initialization that keeps the restarted reflected Halpern iterates uniformly away from zero, with a bound analogous to \cite[Lemma~5.1]{o2021operator}. Recall that $\mathcal Z^\star=\operatorname{Fix}(T)$.

\begin{prop}
\label{prop:hsd-initialization}
Fix $\mu,\eta,\sigma>0$ satisfying \eqref{eq:hsd-step-condition}. Let $x^\circ\in\mathcal X$ and $s^\circ\in\mathcal S$, and let $\{z^{n,t}\}$ be generated by \eqref{eq:restarted-reflected-halpern}--\eqref{eq:restart-point} from
\begin{equation}
    z^{0,0}:=P^{-1}h,
    \qquad
    h:=\bigl(1+c^{\top}x^\circ,\ -c,\ 0,\ s^\circ-Ax^\circ\bigr).
    \label{eq:hsd-anchor}
\end{equation}
Then, for all $n,t$ and every $z^\star=(\tau^\star,x^\star,\alpha^\star,y^\star)\in\mathcal Z^\star\setminus\{0\}$,
\begin{equation}
    \|z^{n,t}\|_P
    \geq
    \frac{\langle z^{0,0},z^\star\rangle_P}{\|z^\star\|_P}
    \geq
    \frac{\tau^\star+\kappa^\star}{\|z^\star\|_P},
    \label{eq:zero-avoidance-bound}
\end{equation}
where $\kappa^\star$ is given by \eqref{eq:hsd-recovery}. In particular,
\begin{equation}
    \inf_{n,t}\|z^{n,t}\|_P>0.
    \label{eq:zero-avoidance-constant}
\end{equation}
\end{prop}

\begin{proof}
Fix $z^\star=(\tau^\star,x^\star,\alpha^\star,y^\star)\in\mathcal Z^\star\setminus\{0\}$, and write $r^\star:=c\tau^\star-A^{\top}y^\star$. The dual-cone constraints in \eqref{eq:four-cone-system} ensure that the support-function values below are finite.

Since $Pz^{0,0}=h$, we have $\langle z^{0,0},z^\star\rangle_P=h^{\top}z^\star$. Expanding and using \eqref{eq:p-func} (with $x^\circ\in\mathcal X$ and $s^\circ\in\mathcal S$) gives
\[
\begin{aligned}
    \langle z^{0,0},z^\star\rangle_P
    &=h^{\top}z^\star\\
    &=\tau^\star-c^{\top}x^\star+(r^\star)^{\top}x^\circ+(y^\star)^{\top}s^\circ\\
    &\geq\tau^\star-c^{\top}x^\star-p(-r^\star;\ell_v,u_v)-p(-y^\star;\ell_c,u_c)\\
    &=\tau^\star+\kappa^\star,
\end{aligned}
\]
where the last equality is \eqref{eq:hsd-recovery}.

We next show that the iteration preserves this inner-product lower bound. Let $z$ be arbitrary, set $u:=T(z)$, and take $v=u+z^\star$ in \eqref{eq:hsd-pdhg-proximal}. Since $u,z^\star\in\mathcal K$ and $\mathcal K$ is a convex cone, we have $v\in\mathcal K$, and thus
\[
    \langle u-z,z^\star\rangle_P
    \geq-\langle Qu,z^\star\rangle
    =\langle u,Qz^\star\rangle
    \geq0.
\]
Here we used $Q^{\top}=-Q$ for the equality, and $u\in\mathcal K$, $Qz^\star\in\mathcal K^*$ for the last inequality. Therefore,
\[
    \langle T(z),z^\star\rangle_P\geq\langle z,z^\star\rangle_P,
    \qquad
    \langle (2T-I)(z),z^\star\rangle_P\geq\langle z,z^\star\rangle_P.
\]

Within outer iteration $n$, applying this to \eqref{eq:restarted-reflected-halpern} yields
\[
    \langle z^{n,t+1},z^\star\rangle_P
    \geq\frac{t+1}{t+2}\langle z^{n,t},z^\star\rangle_P
       +\frac{1}{t+2}\langle z^{n,0},z^\star\rangle_P.
\]
Induction on $t$ gives $\langle z^{n,t},z^\star\rangle_P\geq\langle z^{n,0},z^\star\rangle_P$. At a restart, \eqref{eq:restart-point} gives
\[
    \langle z^{n+1,0},z^\star\rangle_P
    =\langle T(z^{n,t_n}),z^\star\rangle_P
    \geq\langle z^{n,t_n},z^\star\rangle_P
    \geq\langle z^{n,0},z^\star\rangle_P.
\]
Induction on $n$ therefore yields
\begin{equation}
    \langle z^{n,t},z^\star\rangle_P
    \geq\langle z^{0,0},z^\star\rangle_P
    \qquad\text{for all }n,t.
    \label{eq:halpern-inner-product}
\end{equation}

Combining \eqref{eq:halpern-inner-product} with the Cauchy--Schwarz inequality gives
\[
    \|z^{n,t}\|_P\,\|z^\star\|_P\geq\langle z^{n,t},z^\star\rangle_P\geq\langle z^{0,0},z^\star\rangle_P.
\]
Together with the initial inner-product bound, this proves \eqref{eq:zero-avoidance-bound}.

Finally, \autoref{cor:informative-hsd-solution} guarantees a choice of $z^\star$ with $\tau^\star+\kappa^\star>0$. For this fixed solution, $(\tau^\star+\kappa^\star)/\|z^\star\|_P$ is positive and independent of $n$ and $t$, proving \eqref{eq:zero-avoidance-constant}.
\qed
\end{proof}
The proposition excludes convergence to the origin under its stated assumptions. This does not by itself guarantee recovery of an optimal solution or an infeasibility certificate, since a nonzero HSD solution with $\tau=\kappa=0$ may be inconclusive, as discussed in \autoref{subsec:solu}.

The points $x^\circ$ and $s^\circ$ in \eqref{eq:hsd-anchor} can be chosen independently from their respective boxes and need not satisfy $Ax^\circ=s^\circ$. For example, one may take $x^\circ=\Pi_{\mathcal X}(0)$ and $s^\circ=\Pi_{\mathcal S}(0)$.

\subsection{Convergence of the restarted reflected Halpern iteration}
\label{subsec:convergence}

By \eqref{eq:hsd-pdhg-inclusion}, the generalized PDHG operator $T$ is a preconditioned proximal-point map for $Q+N_{\mathcal K}$. This extends the proximal-point interpretation of PDHG~\cite{he2012convergence,lu2025overview} to the HSD monotone inclusion \eqref{eq:hsd-monotone-inclusion}~\cite{o2021operator}. We follow the restarted Halpern analysis of \cite{lu2024restarted}.

Throughout this subsection, the step sizes are fixed and satisfy \eqref{eq:hsd-step-condition}, and $T$ is evaluated exactly. As established in \autoref{subsec:op}, $T$ is firmly nonexpansive in the $P$-norm and $\operatorname{Fix}(T)=\mathcal Z^\star$. The HSD solution set $\mathcal Z^\star$ is nonempty regardless of the status of the original LP. Within each restart cycle, the reflected Halpern iterates satisfy
\begin{equation}
    \|z^{n,t}-T(z^{n,t})\|_P
    \leq
    \frac{\operatorname{dist}_P(z^{n,0},\mathcal Z^\star)}{t+1},
    \qquad 0\leq t\leq t_n.
    \label{eq:hsd-halpern-rate}
\end{equation}
This estimate follows from the Halpern residual bound for nonexpansive operators~\cite[Theorem~2.1]{lieder2021convergence}, using the reflection argument in \cite[Corollary~1]{lu2024restarted}; see also \cite{park2022exact}.

Since $\mathcal K$ and $\mathcal K^*$ are polyhedral, Hoffman's error bound~\cite{hoffman1952approximate} applied to \eqref{eq:box-hsd} gives a constant $C_P>0$ such that the first inequality below holds for every $u\in\mathcal K$. For $u=T(z)$, \eqref{eq:hsd-pdhg-inclusion} gives $u\in\mathcal K$ and $Qu+P(u-z)\in\mathcal K^*$. Using $\|P(u-z)\|_{P^{-1}}=\|u-z\|_P$, we obtain
\[
    \operatorname{dist}_P(u,\mathcal Z^\star)
    \leq
    C_P\operatorname{dist}_{P^{-1}}(Qu,\mathcal K^*)
    \leq C_P\|z-u\|_P.
\]
The triangle inequality then yields
\begin{equation}
    \gamma_P\operatorname{dist}_P(z,\mathcal Z^\star)
    \leq \|z-T(z)\|_P,
    \qquad z\in\mathbb R^{n+m+2},
    \label{eq:hsd-fp-sharpness}
\end{equation}
where $\gamma_P:=(1+C_P)^{-1}>0$.

Combining \eqref{eq:hsd-halpern-rate} and \eqref{eq:hsd-fp-sharpness} with the nonexpansiveness of $T$, we obtain the following convergence guarantee for fixed-frequency restarting.

\begin{thm}[Fixed-frequency restart]
\label{thm:hsd-convergence}
Under the conditions above, consider fixed-frequency restarting with $t_n=\lceil e/\gamma_P\rceil$ and $z^{n+1,0}=T(z^{n,t_n})$ for all $n\geq0$. Then
\begin{equation}
    \operatorname{dist}_P(z^{n,0},\mathcal Z^\star)
    \leq
    e^{-n}\operatorname{dist}_P(z^{0,0},\mathcal Z^\star),
    \qquad n\geq0.
    \label{eq:hsd-linear-rate}
\end{equation}
Moreover, there exists $z^\infty\in\mathcal Z^\star$ such that
\begin{equation}
    \|z^{n,t}-z^\infty\|_P
    \leq
    2e^{-n}\operatorname{dist}_P(z^{0,0},\mathcal Z^\star),
    \qquad n\geq0,\quad 0\leq t\leq t_n.
    \label{eq:hsd-nonzero-linear-rate}
\end{equation}
If $z^{0,0}$ is chosen as in \autoref{prop:hsd-initialization}, then $z^\infty\neq0$.
\end{thm}

\autoref{thm:hsd-convergence} extends the fixed-frequency linear convergence guarantee for restarted reflected Halpern PDHG to the generalized PDHG operator for the HSD embedding. The key ingredients are the firm nonexpansiveness of the operator in the $P$-norm and the global sharpness bound \eqref{eq:hsd-fp-sharpness}. These same ingredients allow the adaptive-restart analysis in~\cite[Theorem~2]{lu2024restarted} to carry over to the present setting with fixed step sizes and exact operator evaluations. We omit the analogous proof to avoid repetition.
The accelerated two-stage analysis in~\cite[Theorems~3 and~4]{lu2024restarted} should also extend by adapting the identification and local convergence arguments to the polyhedral geometry of the HSD system, although the resulting analysis is more involved.

\section{Practical Algorithmic Enhancements}
\label{sec:enhance}

In this section, we describe the practical enhancements used in \HPDLP{}: preprocessing and scaling of the LP data, the step-size parameterization, and the practical initialization (\autoref{subsec:precond}), adaptive restarting and homogeneous rescaling of the restart point (\autoref{subsec:restrt}), and the adaptive primal weight (\autoref{subsec:pid}). With fixed step sizes and exact operator evaluations, \autoref{prop:hsd-initialization} permits arbitrary restart times under the basic restart rule, whereas the estimates in \autoref{thm:hsd-convergence} use the prescribed fixed restart frequency. The practical initialization is not covered by \autoref{prop:hsd-initialization}, and the combined adaptive scheme is not covered by \autoref{thm:hsd-convergence}.

\subsection{Preconditioning, step sizes, and initialization}
\label{subsec:precond}

\HPDLP{} optionally applies the presolver \PSLP{}~\cite{cederberg2026presolving} before the iterative solve. In our experiments, \PSLP{} is applied as a common preprocessing step for all solvers, and results are reported both with and without it (\autoref{subsec:setup}).

\paragraph{Preconditioning}
The problem passed to the iterative solve is then rescaled by the four-stage procedure of \cuPDLPx{}~\cite{lu2025cupdlpx}: $12$ iterations of geometric-mean scaling~\cite{tomlin2009scaling}, $10$ iterations of $\ell_\infty$ Ruiz equilibration~\cite{ruiz2001scaling}, Pock--Chambolle scaling~\cite{pock2011diagonal} with $\alpha=1$, and finally a normalization of the constraint and variable bounds by $\frac{1}{1+\|(\ell_c,u_c)\|_2}$ and of the objective vector by $\frac{1}{1+\|c\|_2}$. Geometric-mean scaling is a classical technique for simplex methods and is available in \Gurobi{}~\cite{gurobi2026}. The combination of Ruiz and Pock--Chambolle scaling was introduced by \PDLP{}~\cite{applegate2026pdlp}. The final normalization is that of \HPRLP{}~\cite{chen2026hpr}. The objective and all finite bounds are transformed consistently under the four diagonal scalings. The algorithm operates on the resulting scaled problem. Inverse scaling, followed by postsolve when \PSLP{} is used, recovers a solution or certificate for the original LP.

\paragraph{Step-size parameterization}
The step-size condition \eqref{eq:hsd-step-condition} has two parts: $\mu\sigma<1$ for the unit coupling between $\tau$ and $\alpha$, and $\eta\sigma\|A\|_2^2<1$ for the coupling between $x$ and $y$ through $A$. We use a common step size $\sigma$ for $\alpha$ and $y$ so that their joint update is a Euclidean projection onto $\mathcal K_c^*$, computed through the Moreau decomposition \eqref{eq:box-cone-moreau}. For $x$ and $y$, we follow \PDLP{}~\cite{applegate2026pdlp} and parameterize $\eta$ and $\sigma$ by a base step size $\lambda>0$ and a primal weight $\omega>0$, which gives $\eta\sigma=\lambda^2$. For $\tau$, we choose $\mu$ so that $\mu\sigma=\nu$ for a fixed $\nu\in(0,1)$. Altogether,
\begin{equation}
    \mu:=\frac{\nu}{\omega\lambda},
    \qquad
    \eta:=\frac{\lambda}{\omega},
    \qquad
    \sigma:=\omega\lambda.
\label{eq:step-size-parameterization}
\end{equation}
The primal weight scales the step sizes of the primal block $(\tau,x)$ by $1/\omega$ and that of the dual block $(\alpha,y)$ by $\omega$. Both products are independent of $\omega$, so \eqref{eq:hsd-step-condition} holds whenever $\lambda\|A\|_2<1$. We set $\nu=0.998$ and $\lambda=\nu/\|A\|_2$, with $\|A\|_2$ estimated by power iteration. We initialize $\omega_0=1$ and keep $\lambda$ and $\nu$ fixed. The primal weight is updated at restarts as described in \autoref{subsec:pid}, and we write $T_n$ and $P_n$ for the operator and the metric used throughout outer iteration $n$.

\paragraph{Practical initialization}
In \HPDLP{}, we use the simpler anchor $z^{0,0}:=(1,0,1,0)$. In our experiments, this anchor required fewer iterations than \eqref{eq:hsd-anchor}, and we observed no convergence to the origin.
\subsection{Adaptive restarting and homogeneous rescaling}
\label{subsec:restrt}

We determine restart times from the fixed-point residual and, when needed, rescale the resulting PDHG point before using it as the next restart anchor.

\paragraph{Adaptive restarting}
We use the fixed-point residual as the restart metric, following~\cite{lu2024restarted}, and apply the three-condition strategy of \cuPDLPx{}~\cite{lu2025cupdlpx}. Define
\begin{equation}
    E^{n,t}
    :=\bigl\|z^{n,t}-T_n(z^{n,t})\bigr\|_{P_n}.
\label{eq:restart-residual}
\end{equation}
Since $T_n(z^{n,t})$ is already computed in the Halpern step, $E^{n,t}$ requires no additional evaluation of $T_n$. Let $K^{n,t}$ denote the total number of inner iterations performed up to iterate $(n,t)$. A restart is triggered if any of the following holds:
\begin{equation}
\begin{aligned}
    E^{n,t}
    &\leq \beta_{\mathrm{su}}E^{n,0}
    &&\text{(sufficient decay)},\\
    E^{n,t}
    &\leq \beta_{\mathrm{ne}}E^{n,0}
    \quad\text{and}\quad E^{n,t}>E^{n,t-200}
    &&\text{(necessary decay without local progress)},\\
    t
    &\geq \beta_{\mathrm{a}}K^{n,t}
    &&\text{(long inner loop)}.
\end{aligned}
\label{eq:adaptive-restart-conditions}
\end{equation}
We use $\beta_{\mathrm{su}}=0.2$, $\beta_{\mathrm{ne}}=0.5$, and $\beta_{\mathrm{a}}=0.36$. The conditions are evaluated every $200$ inner iterations, so the second condition compares the residual with its value at the previous evaluation. If a restart is triggered at inner iteration $t_n$, we write
\begin{equation}
    \bar z_n:=T_n(z^{n,t_n})=(\bar\tau_n,\bar x_n,\bar\alpha_n,\bar y_n)
\label{eq:restart-pdhg-point}
\end{equation}
for the PDHG point at the restart.

\paragraph{Homogeneous rescaling}
The HSD formulation does not fix the scale of a solution: every positive multiple of a solution is again a solution. We rescale restart points to keep their magnitude within a suitable numerical range. Before rescaling, define
\begin{equation}
    h_n:=\max\bigl\{\|\bar x_n\|_2,\ \|\bar y_n\|_2,\ |\bar\alpha_n|\bigr\},
    \qquad
    g_n:=\max\{h_n,\ \bar\tau_n\},
\label{eq:restart-point-magnitude}
\end{equation}
and set the anchor of the next outer iteration to
\begin{equation}
    z^{n+1,0}:=\gamma_n\bar z_n,
    \qquad
    \gamma_n:=
    \begin{cases}
        \frac{1}{g_n}, & g_n\notin[10^{-12},10^{12}],\\
        \frac{1}{\bar\tau_n}, & \text{otherwise, if } \bar\tau_n>10^{-8}h_n \text{ and } \bar\tau_n\notin[0.1,10],\\
        1, & \text{otherwise}.
    \end{cases}
\label{eq:rescaled-restart-point}
\end{equation}
The first case rescales points with extreme magnitudes to reduce the risk of overflow and underflow. Otherwise, we normalize $\bar\tau_n$ to one only when it lies outside $[0.1,10]$ and is not negligible relative to the other components. This makes the anchor's $(x,y)$-components equal to the recovered primal-dual point while avoiding excessive amplification when $\bar\tau_n$ is small near the certificate branch. In all other cases, we leave the restart point unchanged.

Positive rescaling preserves the recovered primal-dual point $(x/\tau,y/\tau)$ when $\tau>0$ and the certificate directions. For a fixed operator $T_n$, positive homogeneity implies that rescaling the anchor scales the inner iterates and their residuals by the same factor, preserving the relative restart comparisons.

\subsection{Adaptive primal weight}
\label{subsec:pid}

The primal weight is intended to balance the weighted primal and dual distances to a solution. Since the solution is unknown, we use the movement during the completed outer iteration as a proxy, as in \PDLP{} and \cuPDLPx{}~\cite{applegate2026pdlp,lu2025cupdlpx}.

\paragraph{Target primal weight}
At the end of outer iteration $n$, before homogeneous rescaling, define
\begin{equation}
    \Delta z_n
    :=\bar z_n-z^{n,0}
    =(\Delta\tau_n,\Delta x_n,\Delta\alpha_n,\Delta y_n),
\label{eq:restart-movement}
\end{equation}
and the movements of the primal block $(\tau,x)$ and the dual block $(\alpha,y)$,
\begin{equation}
    (\delta^{\mathrm p}_n)^2
    :=\|\Delta x_n\|_2^2+\frac{\lambda^2}{\nu}\Delta\tau_n^2,
    \qquad
    (\delta^{\mathrm d}_n)^2
    :=\|\Delta y_n\|_2^2+\Delta\alpha_n^2.
\label{eq:primal-dual-movement}
\end{equation}
The factor $\frac{\lambda^2}{\nu}=\frac{\eta}{\mu}$ accounts for the different step sizes of $\tau$ and $x$. Writing $P(\omega)$ for the metric \eqref{eq:hsd-pdhg-metric} under the parameterization \eqref{eq:step-size-parameterization},
\begin{equation}
    \|\Delta z_n\|_{P(\omega)}^2
    =
    \frac{\omega}{\lambda}(\delta^{\mathrm p}_n)^2
    +\frac{1}{\omega\lambda}(\delta^{\mathrm d}_n)^2
    -2\Delta\tau_n\Delta\alpha_n
    -2\langle A\Delta x_n,\Delta y_n\rangle.
\label{eq:primal-weight-movement}
\end{equation}
The last two terms do not depend on $\omega$. When $\delta^{\mathrm p}_n,\delta^{\mathrm d}_n>0$, minimizing the $\omega$-dependent terms over $\omega>0$ gives the target
\begin{equation}
    \omega_n^\star:=\frac{\delta^{\mathrm d}_n}{\delta^{\mathrm p}_n}.
\label{eq:target-primal-weight}
\end{equation}

\paragraph{Gating in the certificate regime}
Near a certificate solution, the movements of the primal and dual blocks may differ substantially, making their ratio an unreliable target for balancing the step sizes. We therefore suspend weight adaptation when the restart point indicates the certificate regime. Let $\bar r_n:=\Pi_{\mathcal R}(c\bar\tau_n-A^{\top}\bar y_n)$ be the projected reduced cost at $\bar z_n$, and let
\begin{equation}
    \kappa_n
    :=-c^{\top}\bar x_n
      -p(-\bar y_n;\ell_c,u_c)
      -p(-\bar r_n;\ell_v,u_v)
\label{eq:certificate-scalar-restart}
\end{equation}
be the corresponding value of $\kappa$ in \eqref{eq:hsd-recovery}. Both support-function values are finite, the first because $(\bar\alpha_n,\bar y_n)\in\mathcal K_c^*$ and the second because $\bar r_n\in\mathcal R$. Define
\begin{equation}
    \rho_n
    :=
    \begin{cases}
        \dfrac{(\kappa_n)^+}{\bar\tau_n+(\kappa_n)^+},
        & \bar\tau_n+(\kappa_n)^+>0,\\[2ex]
        0, & \text{otherwise}.
    \end{cases}
\label{eq:pid-regime-indicator}
\end{equation}
At an HSD solution with $\tau+\kappa>0$, $\rho_n$ equals $0$ on the optimal branch and $1$ on the certificate branch. We suspend the update when $\rho_n>0.9$.

\paragraph{Primal-weight update}
Following the PID primal-weight update of \cuPDLPx{}~\cite{lu2025cupdlpx}, we control the logarithmic imbalance
\begin{equation}
    e_n
    :=\log\frac{\delta^{\mathrm d}_n}{\omega_n\delta^{\mathrm p}_n}
    =\log\frac{\omega_n^\star}{\omega_n}.
\label{eq:primal-weight-error}
\end{equation}
We use the proportional and integral terms, that is, the \cuPDLPx{} update with the derivative term switched off, and discount the integral state exponentially. With $s_{-1}:=0$,
\begin{equation}
    s_n:=\zeta s_{n-1}+e_n,
    \qquad
    \omega_{n+1}
    :=\omega_n\exp\bigl(k_{\mathrm P}e_n+k_{\mathrm I}s_n\bigr),
\label{eq:primal-weight-pid}
\end{equation}
where $k_{\mathrm P},k_{\mathrm I}\geq0$ are the proportional and integral gains and $\zeta\in[0,1)$ is the discount factor. With $k_{\mathrm I}=0$, the update reduces to the exponential smoothing $\omega_{n+1}=\omega_n^{1-k_{\mathrm P}}(\omega_n^\star)^{k_{\mathrm P}}$ of \PDLP{}, which uses $k_{\mathrm P}=0.5$. We use $k_{\mathrm P}=0.97$, $k_{\mathrm I}=0.05$, and $\zeta=0.3$.

If the update is suspended, we set $\omega_{n+1}:=\omega_n$ and $s_n:=0$; otherwise, if $\delta^{\mathrm p}_n=0$ or $\delta^{\mathrm d}_n=0$, we set $\omega_{n+1}:=\omega_n$ and $s_n:=s_{n-1}$. Whenever $\omega$ changes, we recompute the step sizes from \eqref{eq:step-size-parameterization} and hold them fixed throughout the next outer iteration.

\section{A Benchmark Collection of Infeasible and Unbounded LPs}
\label{sec:bench}

The \netlib{} repository~\cite{gay1985electronic} provides a standard collection of LPs with finite optimal values and a separate collection of 29 infeasible instances. The infeasible instances are relatively small, limiting their usefulness for evaluating solver scalability, and the repository has no corresponding collection of unbounded LPs. This motivates a larger benchmark for evaluating infeasibility and unboundedness detection.

We introduce \MIPLIBInfeas{}, a collection of 334 LPs derived from \MIPLIB{} 2017~\cite{gleixner2021miplib}. It contains 229 primal-infeasible instances obtained from branch-and-bound node relaxations and 105 feasible but unbounded instances obtained from root relaxations. \autoref{tab:collection-comparison} compares the size and scale of the collection with those of the \netlib{} infeasible collection.

\begin{table}[!htb]
\centering
\caption{Comparison with the \netlib{} infeasible collection. Rows, Cols, and Nnzs denote the numbers of constraints, variables, and nonzeros; each maximum is computed separately over the instances in the corresponding table row.}
\label{tab:collection-comparison}
\footnotesize
\begin{adjustbox}{center}
\begin{tabular}{llrrrr}
\toprule
& & & \multicolumn{3}{c}{Maximum} \\
\cmidrule(lr){4-6}
Collection & Status & \# & Rows & Cols & Nnzs \\
\midrule
\netlib{} & infeasible & 29 & 3{,}793 & 10{,}733 & 97{,}257 \\
\MIPLIBInfeas{} & infeasible & 229 & 19{,}912{,}111 & 20{,}682{,}487 & 47{,}513{,}027 \\
& unbounded & 105 & 19{,}335{,}696 & 550{,}539 & 58{,}613{,}680 \\
\bottomrule
\end{tabular}
\end{adjustbox}
\end{table}

The collection is available at
\begin{center}
    \url{https://mit-lu-lab.github.io/MIPLIB2017-Infeas}
\end{center}
The website provides individual and bulk downloads in compressed MPS format, together with metadata identifying the source model and recording the modifications used to construct each instance. Searchable and filterable tables report instance dimensions and nonzero counts, while individual instance pages display sparsity patterns and construction details.

\subsection{Construction from MIPLIB}
\label{subsec:constr}

We use the 383 \MIPLIB{} 2017 models whose root-node LP relaxations form the standard \MIPLIB{} relaxations benchmark used in \cuPDLP{} and \cuPDLPx{}~\cite{lu2025cupdlp,lu2025cupdlpx}. The \cuPDLP{} benchmark excludes instances tagged as numerically unstable, infeasible, or having indicator constraints; instances with at most $10^5$ nonzeros; and relaxations for which the zero vector is an optimal solution. Following \cuPDLPx{}, we classify the sources as small, medium, or large according to the number of nonzeros in their constraint matrices: fewer than $10^6$, between $10^6$ and $10^7$, or more than $10^7$. The three groups contain 269, 94, and 20 models.

Using \cplex{}~22.1.2~\cite{manual1987ibm}, we generate candidates separately for infeasibility and unboundedness. We discard generated candidates whose status is detected by the presolver of either \HiGHS{}~1.14.0 or \Gurobi{}~12.0.3. For each source and status class, we retain up to five, three, or one candidates for small, medium, or large sources, respectively, within a time limit of 10{,}800 seconds. We then sample one of the retained candidates using a fixed random seed; if none are retained, no instance is selected. Each source thus contributes at most one infeasible LP and at most one unbounded LP. All released instances are continuous LPs, with integrality restrictions removed. Each selected instance is reloaded from its distributed MPS file and verified with \Gurobi{}, requiring a definitive infeasible or unbounded status.

\paragraph{Infeasible instances}
We obtain infeasible LPs from branch-and-bound searches on the original MIP models. During the branch-and-bound search, cuts, presolve, probing, and primal heuristics are disabled, so each node LP is obtained from the root relaxation by tightening variable bounds through branching. Since the \cplex{} relaxation callback does not expose infeasible node relaxations directly, we construct sibling LPs from visited nodes. Whenever the branching bound can be identified, we apply its complementary bound to the parent's variable bounds to form the sibling relaxation. We retain sibling LPs that pass the presolve screening described above and are subsequently verified as infeasible. The selected instances correspond to nodes at depths ranging from 2 to 54{,}505, with a median of 149; 133 of the 229 instances have depth at least 100.

\paragraph{Unbounded instances}
The root LP relaxations of \texttt{neos-3229051-yass} and \texttt{neos-3230511-yuna} are already unbounded. We include them unchanged, without candidate generation or presolve screening. For each remaining feasible and bounded root relaxation, we remove constraints and finite variable bounds one at a time, re-solving the LP after each removal until it becomes unbounded. At each step, eligible constraints have nonzero optimal dual multipliers, and eligible bounds are active with nonzero reduced costs. We assign each eligible constraint a weight equal to the absolute value of its dual multiplier times its row nonzero count, and each eligible bound a weight equal to the absolute reduced cost times its column nonzero count. We sample one constraint or bound with probability proportional to its weight and remove it. Once the LP becomes unbounded, we test the removed constraints and bounds one at a time, restoring each if the LP remains unbounded. The resulting LP remains feasible because it is a relaxation of the original feasible LP. Across the 105 unbounded instances, the median numbers of removed constraints and variable bounds are 1 and 5.

\subsection{Scale and source coverage}
\label{subsec:scale}

The collection covers 262 of the 383 source models, with 72 contributing both an infeasible and an unbounded instance. Infeasible instances were obtained from 171 of the 269 small, 52 of the 94 medium, and 6 of the 20 large sources. The corresponding counts for unbounded instances are 93, 9, and 3, respectively. For both classes, coverage is lower among medium and large sources than among small sources.

\begin{figure}[!htb]
\centering
\includegraphics[width=0.6\linewidth]{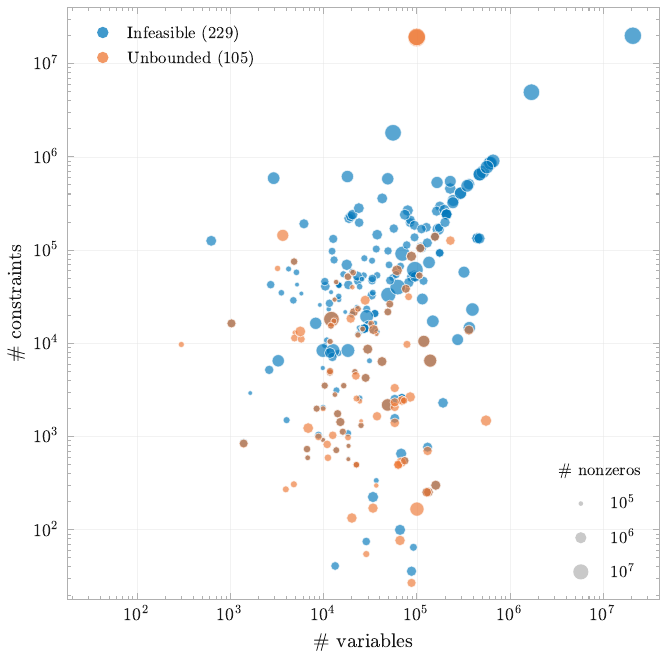}
\caption{Scale of the benchmark instances. Each marker represents one LP; color indicates its status class and marker size its number of nonzeros. Both axes are logarithmic.}
\label{fig:infeas-scatter}
\end{figure}

Across both classes, the instances range from 297 to 20{,}682{,}487 variables, from 27 to 19{,}912{,}111 constraints, and from 100{,}933 to 58{,}613{,}680 nonzeros; 64 instances have at least $10^6$ nonzeros. \autoref{fig:infeas-scatter} shows the joint distribution of variables and constraints, with marker size indicating the number of nonzeros.

\section{Numerical Experiments}
\label{sec:exp}

We first evaluate infeasibility and unboundedness detection, then compare performance on feasible LP benchmarks. We also examine iteration counts and computational costs, including the box-cone subproblems. The numerical comparisons aim to identify the strengths and limitations of \HPDLP{}, rather than to establish the superiority of \HPDLP{} over existing LP approaches.

\subsection{Experimental setup}
\label{subsec:setup}

We compare \HPDLP{} with CPU and GPU solvers, both with and without preprocessing by \PSLP{}~0.0.8~\cite{cederberg2026presolving}. All solvers receive the same reduced problems in the presolved experiments, and all solver-specific presolvers are disabled in both settings.

\paragraph{Benchmark sets}
For infeasibility and unboundedness detection, we use all 334 instances in \MIPLIBInfeas{} (\autoref{sec:bench}) and all 29 instances in the \netlib{} infeasible collection~\cite{gay1985electronic}. We also evaluate performance on two feasible LP benchmark sets, namely, the 383 LP relaxations of \MIPLIB{} 2017 models~\cite{gleixner2021miplib,lu2025cupdlp} and the 49 public instances of Mittelmann's LP benchmark~\cite{mittelmann2009decision,lu2025cupdlpx}. The \MIPLIB{} relaxations use the size groups defined in \autoref{subsec:constr} and include the two originally unbounded LPs identified there.
 
\paragraph{Software and algorithms}
The baseline solvers are \Gurobi{}~12.0.3~\cite{gurobi2026}, \HiGHS{}~1.14.0~\cite{huangfu2018parallelizing}, \MOSEK{}~11.2.2~\cite{aps2022mosek}, \SCS{}~3.2.11~\cite{o2016conic}, the C implementation of \ABIP{}~2.0.0~\cite{deng2022enhanced}, the \PDLP{} implementation in Google \ORTools{}~9.15.6755~\cite{applegate2026pdlp,ortools}, \cuPDLPx{}~0.3.0~\cite{lu2025cupdlpx}, and \HPRLPC{}~0.1.2, the C/GPU implementation of \HPRLP{}~\cite{chen2026hpr}. For \Gurobi{}, we test primal simplex, dual simplex, barrier, and homogeneous barrier, with crossover disabled for both barrier variants. \HiGHS{} uses dual simplex, \MOSEK{} uses its interior-point optimizer, and \SCS{} uses a direct linear system solver. \cuPDLPx{} and \HPRLPC{} are evaluated only on the feasible LP benchmarks because they do not provide infeasibility detection. Unless otherwise specified, solver parameters retain their default values. The iteration limit of \SCS{} is increased to $10^9$ to avoid termination at the iteration limit.
 
\paragraph{Computing environment}
All experiments are conducted on Rocky Linux~8.10. Each run is allocated 128~GB of host memory. CPU runs are allocated 16 AMD EPYC~9474F cores, while GPU runs are allocated one CPU core and one NVIDIA H100 GPU with 80~GB of device memory. \HPDLP{}, \cuPDLPx{}, and \HPRLPC{} are compiled with GCC~12.2.0 and CUDA Toolkit~13.3.1.

\subsection{Termination criteria and evaluation}
\label{subsec:term}

\HPDLP{} checks for primal-dual optimality and certificates of primal or dual infeasibility. In the tests below, $(\tau,x,\alpha,y)$ denotes the PDHG point $T(z^{n,t})$ after undoing the scaling of \autoref{subsec:precond}. After describing these tests, we specify the baseline tolerances, time limits, and performance metrics.

\paragraph{Optimality}
For $\tau>0$, let $q:=c\tau-A^{\top}y$. Let $\pi$ denote the variable-bound multipliers from the primal box-cone projection, with positive entries for lower bounds and negative entries for upper bounds. We consider two reduced-cost candidates:
\[
    r_j^{\mathrm{act}}:=
    \begin{cases}
        q_j, & \pi_jq_j>0,\\
        0,   & \text{otherwise},
    \end{cases}
    \qquad
    r^{\mathrm{bnd}}:=\Pi_{\mathcal R}(q).
\]
The second follows the projected reduced-cost construction used in \cuPDLPC{}~\cite{lu2023cupdlp}. For each candidate $r$, set $(\widehat x,\widehat y,\widehat r):=(x,y,r)/\tau$, with $\widehat x\in\mathcal X$ and $\widehat r\in\mathcal R$ by construction. Optimality is declared if all three conditions below hold for at least one candidate:
\begin{equation}
\begin{aligned}
    \frac{\|A\widehat x-\Pi_{\mathcal S}(A\widehat x)\|_2}
         {1+\|(\ell_c,u_c)\|_2}
    &\leq\varepsilon,\\
    \frac{\|A^{\top}\widehat y+\widehat r-c\|_2}
         {1+\|c\|_2}
    &\leq\varepsilon,\\
    \frac{\bigl|c^{\top}\widehat x+p(-\widehat y;\ell_c,u_c)+p(-\widehat r;\ell_v,u_v)\bigr|}
         {1+|c^{\top}\widehat x|+\bigl|p(-\widehat y;\ell_c,u_c)+p(-\widehat r;\ell_v,u_v)\bigr|}
    &\leq\varepsilon.
\end{aligned}
\label{eq:optimality-tolerance}
\end{equation}
These are the relative KKT criteria of \PDLP{}~\cite{applegate2026pdlp} in the bound form of \cuPDLPx{}~\cite{lu2025cupdlpx}. On the feasible LP benchmarks, we run \HPDLP{}, \cuPDLPx{}, and \HPRLPC{} at the same nominal tolerances, $\varepsilon=10^{-4}$ and $\varepsilon=10^{-8}$, setting both the feasibility and the optimality tolerance of \HPDLP{} and \cuPDLPx{} to $\varepsilon$.

\paragraph{Infeasibility certificates}
Motivated by \autoref{prop:hsd-solution-characterization}, we test $(y,r_{\mathrm{cert}})$, with $r_{\mathrm{cert}}:=\Pi_{\mathcal R}(-A^{\top}y)$, as a candidate certificate of primal infeasibility \eqref{eq:primal-infeas-cert}, and $x$ as a candidate certificate of dual infeasibility \eqref{eq:dual-infeas-cert}. We also form candidates from the displacement of the PDHG point from the current restart anchor $z^{n,0}$:
\begin{equation}
    d_x:=\Pi_{\operatorname{rec}(\mathcal X)}(x-x^{n,0}),
    \qquad
    d_y:=\Pi_{\mathcal Y}(y-y^{n,0}),
    \qquad
    d_r:=\Pi_{\mathcal R}(-A^{\top}d_y).
\label{eq:displacement-candidates}
\end{equation}
These projections ensure that $d_x\in\operatorname{rec}(\mathcal X)$, $d_y\in\mathcal Y$, and $d_r\in\mathcal R$. The displacement construction is motivated by the certificates obtained from PDHG iterates~\cite{applegate2024infeasibility}, but we do not establish a corresponding convergence guarantee for these candidates.
 
The primal and dual infeasibility tests are, respectively,
\begin{equation}
\begin{aligned}
    \frac{\|A^{\top}y+r_{\mathrm{cert}}\|_\infty}
         {-p(-y;\ell_c,u_c)-p(-r_{\mathrm{cert}};\ell_v,u_v)}
    &\leq\varepsilon_{\mathrm{cert}},\\
    \frac{\max\bigl\{
        \|x-\Pi_{\operatorname{rec}(\mathcal X)}(x)\|_\infty,\,
        \|Ax-\Pi_{\operatorname{rec}(\mathcal S)}(Ax)\|_\infty
        \bigr\}}
         {-c^{\top}x}
    &\leq\varepsilon_{\mathrm{cert}}.
\end{aligned}
\label{eq:certificate-tolerance}
\end{equation}
A candidate is accepted only if its denominator is positive and the corresponding inequality holds. Tests of this scale-invariant form are standard for homogeneous embeddings~\cite{andersen2000primal,o2016conic}. Each accepted certificate is normalized so that its denominator equals one. We use $\varepsilon_{\mathrm{cert}}=10^{-9}$ for $(y,r_{\mathrm{cert}})$ and $x$, and apply the same tests to the displacement candidates $(d_y,d_r)$ and $d_x$ with $\varepsilon_{\mathrm{cert}}=10^{-10}$. These are the default certificate tolerances of \HPDLP{} and are independent of the optimality tolerance $\varepsilon$.

\paragraph{Baseline tolerances}
The CPU baselines use their native stopping criteria at default tolerances, following the \PDLP{} benchmark study~\cite{applegate2026pdlp}. The exception is \ABIP{}, whose tolerance is set to $10^{-4}$ instead of the default $10^{-3}$, following its benchmark experiments~\cite{deng2022enhanced}. \autoref{tab:tolerances} summarizes the stopping criteria and tolerance settings. Since these criteria differ across solvers, equal numerical tolerances need not imply equal solution accuracy~\cite{lu2025cupdlp}.
 
\begin{table}[!htb]
\centering
\caption{Stopping criteria and tolerances. Defaults are used except for the optimality tolerances of \ABIP{} and the GPU solvers. ``abs.'' uses unnormalized residuals, while ``rel.'' measures residuals and gaps relative to data and objective magnitudes. ``ratio'' divides a certificate residual by its objective magnitude, as in \eqref{eq:certificate-tolerance}. The GPU solvers use optimality tolerances of $10^{-4}$ and $10^{-8}$; the certificate tolerances of \HPDLP{} remain fixed. \Gurobi{} barrier runs without crossover.}
\label{tab:tolerances}
\footnotesize
\setlength{\tabcolsep}{2.5pt}
\begin{adjustbox}{center}
\begin{threeparttable}
\begin{tabular}{llll}
\toprule
& \multicolumn{2}{c}{Optimality tests} & Infeasibility detection \\
\cmidrule(lr){2-3}\cmidrule(lr){4-4}
Solver & Feasibility & Gap & Criterion \\
\midrule
\Gurobi{} primal simplex & abs. $\leq10^{-6}$ & --- & from final basis \\
\Gurobi{} dual simplex & abs. $\leq10^{-6}$ & --- & from final basis \\
\Gurobi{} barrier & rel. $\leq10^{-8}$ & rel.\tnote{a}\hspace{0.8em}$\leq10^{-8}$ & status only \\
\Gurobi{} homogeneous barrier & rel. $\leq10^{-8}$ & rel.\tnote{a}\hspace{0.8em}$\leq10^{-8}$ & status only \\
\HiGHS{} dual simplex & abs. $\leq10^{-7}$ & --- & from final basis \\
\MOSEK{} interior point & rel.\ $\ell_\infty$ $\leq10^{-8}$ & rel. $\leq10^{-8}$ & $\ell_\infty$ ratio $\leq10^{-10}$ \\
\SCS{} & abs.\,$+$\,rel.\ $\ell_\infty$ $\leq10^{-4}$ & abs.\,$+$\,rel. $\leq10^{-4}$ & $\ell_\infty$ ratio $\leq10^{-7}$ \\
\ABIP{} & rel.\ $\ell_2$ $\leq10^{-4}$ & rel. $\leq10^{-4}$ & $\ell_2$ ratio $\leq10^{-4}$ \\
\ORTools{} \PDLP{} & abs.\,$+$\,rel.\ $\ell_2$ $\leq10^{-6}$ & abs.\,$+$\,rel. $\leq10^{-6}$ & $\ell_\infty$ ratio $\leq10^{-8}$ \\
\midrule
\HPRLPC{} & rel.\ $\ell_2$ $\leq10^{-4}/10^{-8}$ & rel. $\leq10^{-4}/10^{-8}$ & ---  \\
\cuPDLPx{} & rel.\ $\ell_2$ $\leq10^{-4}/10^{-8}$ & rel. $\leq10^{-4}/10^{-8}$ & --- \\
\HPDLP{} & rel.\ $\ell_2$ $\leq10^{-4}/10^{-8}$ & rel. $\leq10^{-4}/10^{-8}$ & $\ell_\infty$ ratio\tnote{b}\hspace{0.8em}$\leq10^{-9}$ \\
\bottomrule
\end{tabular}
\begin{tablenotes}
\footnotesize
\item[a] Also bounds complementarity.
\item[b] $10^{-10}$ for the displacement candidates \eqref{eq:displacement-candidates}.
\end{tablenotes}
\end{threeparttable}
\end{adjustbox}
\end{table}

\paragraph{Time limits}
The solve-time limit is 3{,}600 seconds for \MIPLIBInfeas{}, \netlib{}, and the small and medium \MIPLIB{} relaxations; 18{,}000 seconds for the large \MIPLIB{} relaxations; and 15{,}000 seconds for Mittelmann's instances. Presolve time is excluded from both the reported solve times and these limits. Instances resolved by \PSLP{} are assigned zero solve time.

\paragraph{Performance metrics}
For each benchmark set, we report the outcome counts and the shifted geometric mean of solve times with a shift of 10 seconds,
\[
\operatorname{SGM10}
:=\left(\prod_{i=1}^{N}(t_i+10)\right)^{\frac{1}{N}}-10,
\]
where $N$ is the number of instances. For runs that return an accepted infeasibility status or terminate as optimal under the solver's numerical stopping criteria, we use the actual termination time; all other outcomes are charged the time limit. Thus, SGM10 measures termination time under the prescribed stopping criteria.

On the infeasible and unbounded benchmarks, Detected counts reports consistent with the known instance status, including infeasible-or-unbounded reports and detections by \PSLP{}. Opt.\ counts terminations satisfying the solver's numerical optimality criteria at the prescribed tolerance. On the feasible LP benchmarks, Solved counts both types of accepted termination. Limit denotes runs reaching a time or iteration limit, and Fail denotes all remaining outcomes.

\subsection{Infeasible and unbounded benchmarks}
\label{subsec:infeas}

We evaluate the ability of \HPDLP{} to detect infeasibility and unboundedness on \MIPLIBInfeas{} and the infeasible \netlib{} collection. \autoref{tab:miplib-infeas-results} and \autoref{tab:netlib-infeas-results} report termination outcomes and solve times, with established solvers included as reference points.

\begin{table}[!htb]
\centering
\caption{Results on the 334 instances of \MIPLIBInfeas{}. Outcome categories and SGM10 are defined in \autoref{subsec:term}. Each optimality tolerance of \HPDLP{} is compared separately with the CPU baselines. Within each presolve setting, the highest Detected count and lowest SGM10 in either comparison are bolded; those attaining the best values in the $10^{-8}$ comparison are additionally underlined. Ties at the displayed precision receive the same formatting.}
\label{tab:miplib-infeas-results}
\footnotesize
\setlength{\tabcolsep}{2.5pt}
\begin{adjustbox}{width=0.85\linewidth}
\begin{threeparttable}
\begin{tabular}{lrrrrrrrrrr}
\toprule
& \multicolumn{5}{c}{Without presolve}
& \multicolumn{5}{c}{With \PSLP{}} \\
\cmidrule(lr){2-6}\cmidrule(lr){7-11}
Solver & Detected & Opt. & Limit & Fail & SGM10 (sec) & Detected & Opt. & Limit & Fail & SGM10 (sec) \\
\midrule
\Gurobi{} primal simplex & 323 & 0 & 9 & 2 & 11.52 & \textbf{\underline{334}} & 0 & 0 & 0 & 4.55 \\
\Gurobi{} dual simplex   & \textbf{\underline{332}} & 0 & 0 & 2 & \textbf{\underline{6.05}} & 332 & 0 & 0 & 2 & 2.98 \\
\Gurobi{} barrier        & 187 & 4 & 1 & 142 & 153.63 & 253 & 1 & 0 & 80 & 39.65 \\
\Gurobi{} homogeneous barrier & 257 & 1 & 0 & 76\tnote{b} & 42.84 & 288 & 1 & 0 & 45 & 15.59 \\
\HiGHS{}                 & 330 & 0 & 4 & 0 & 14.47 & 330 & 0 & 1 & 3 & 6.23 \\
\MOSEK{}                 & 327 & 3 & 2 & 2\tnote{a} & 6.22 & 333 & 0 & 0 & 1\tnote{a} & 3.16 \\
\SCS{}                   & 187 & 47 & 94 & 6\tnote{a} & 264.62 & 246 & 33 & 51 & 4\tnote{a} & 84.00 \\
\ABIP{} ($10^{-4}$)      & 82 & 58 & 102 & 92 & 723.84 & 116 & 60 & 45 & 113 & 375.20 \\
\ORTools{} \PDLP{}       & 243 & 4 & 87 & 0 & 305.34 & 277 & 7 & 50 & 0 & 103.99 \\
\midrule
\HPDLP{} ($10^{-4}$)   & 292 & 41 & 1 & 0 & \textbf{3.79} & 305 & 29 & 0 & 0 & \textbf{1.35} \\
\midrule
\HPDLP{} ($10^{-8}$)   & 326 & 3 & 5 & 0 & 6.17 & \textbf{\underline{334}} & 0 & 0 & 0 & \textbf{\underline{2.10}} \\
\bottomrule
\end{tabular}
\begin{tablenotes}
\footnotesize
\item[a] Includes inconsistent infeasibility reports.
\item[b] Includes 1 inconsistent infeasibility report and 75 numerical issues or out-of-memory terminations.
\end{tablenotes}
\end{threeparttable}
\end{adjustbox}
\end{table}

\begin{table}[!htb]
\centering
\caption{Results on the 29 infeasible \netlib{} instances. Column definitions and highlighting conventions follow \autoref{tab:miplib-infeas-results}. All Fail outcomes are numerical statuses or crashes.}
\label{tab:netlib-infeas-results}
\footnotesize
\setlength{\tabcolsep}{2.5pt}
\begin{adjustbox}{width=0.85\linewidth}
\begin{tabular}{lrrrrrrrrrr}
\toprule
& \multicolumn{5}{c}{Without presolve}
& \multicolumn{5}{c}{With \PSLP{}} \\
\cmidrule(lr){2-6}\cmidrule(lr){7-11}
Solver & Detected & Opt. & Limit & Fail & SGM10 (sec) & Detected & Opt. & Limit & Fail & SGM10 (sec) \\
\midrule
\Gurobi{} primal simplex & \textbf{\underline{29}} & 0 & 0 & 0 & \textbf{\underline{0.01}} & \textbf{\underline{29}} & 0 & 0 & 0 & 0.01 \\
\Gurobi{} dual simplex   & \textbf{\underline{29}} & 0 & 0 & 0 & 0.04 & \textbf{\underline{29}} & 0 & 0 & 0 & \textbf{\underline{0.00}} \\
\Gurobi{} barrier        & 23 & 0 & 0 & 6 & 24.01 & 27 & 1 & 0 & 1 & 2.28 \\
\Gurobi{} homogeneous barrier & 27 & 0 & 0 & 2 & 5.05 & \textbf{\underline{29}} & 0 & 0 & 0 & 0.01 \\
\HiGHS{}                 & 27 & 0 & 0 & 2 & 5.07 & \textbf{\underline{29}} & 0 & 0 & 0 & \textbf{\underline{0.00}} \\
\MOSEK{}                 & \textbf{\underline{29}} & 0 & 0 & 0 & 0.04 & \textbf{\underline{29}} & 0 & 0 & 0 & \textbf{\underline{0.00}} \\
\SCS{}                   & 20 & 9 & 0 & 0 & 3.93 & 26 & 3 & 0 & 0 & 0.61 \\
\ABIP{} ($10^{-4}$)      & 17 & 4 & 1 & 7 & 45.29 & 23 & 3 & 0 & 3 & 9.69 \\
\ORTools{} \PDLP{}       & 24 & 0 & 5 & 0 & 74.55 & 27 & 0 & 2 & 0 & 21.97 \\
\midrule
\HPDLP{} ($10^{-4}$)   & 28 & 1 & 0 & 0 & 1.03 & 28 & 1 & 0 & 0 & 0.25 \\
\midrule
\HPDLP{} ($10^{-8}$)   & \textbf{\underline{29}} & 0 & 0 & 0 & 2.01 & \textbf{\underline{29}} & 0 & 0 & 0 & 1.21 \\
\bottomrule
\end{tabular}
\end{adjustbox}
\end{table}

\begin{figure}[!htb]
    \centering
    \includegraphics[width=0.9\linewidth]{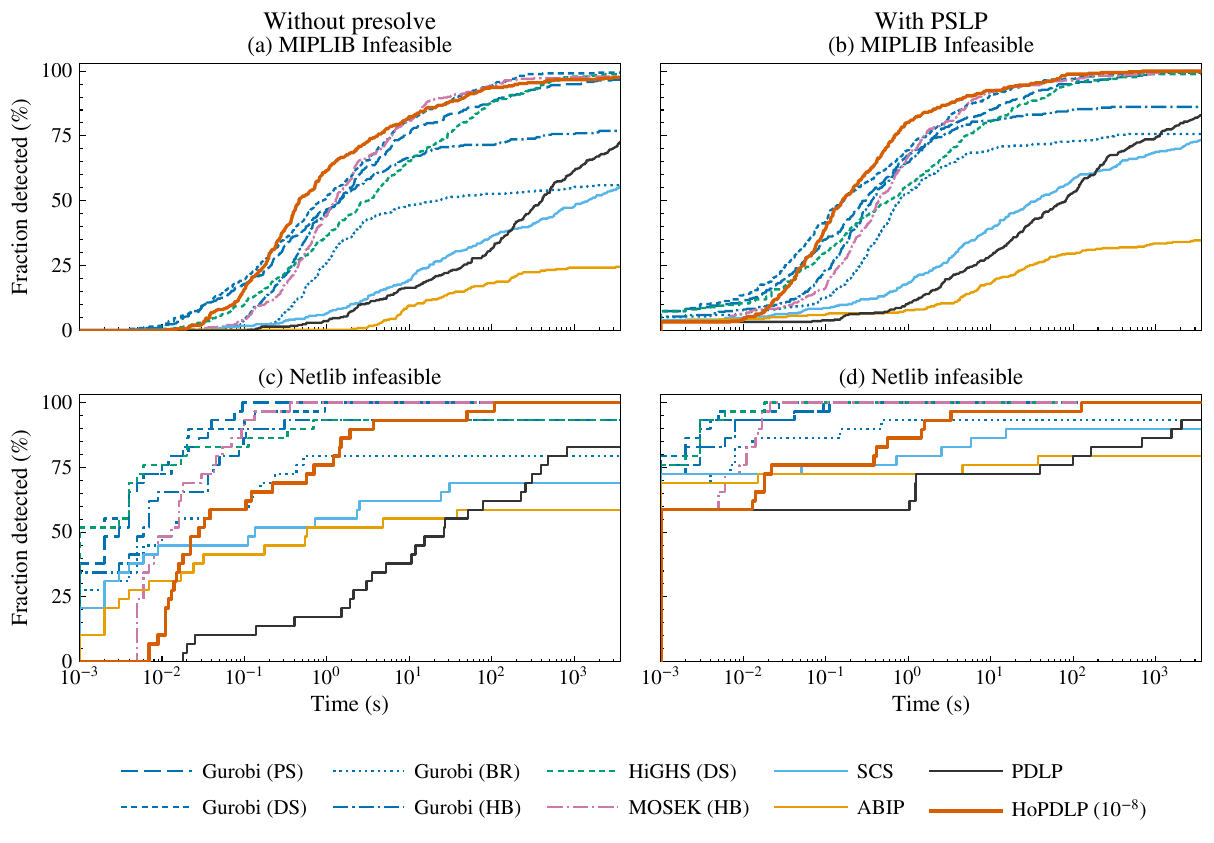}
    \caption{Cumulative detection rates without presolve and with \PSLP{}. \HPDLP{} uses $\varepsilon=10^{-8}$. Fractions are relative to all instances in each collection, counting detections within 3600 seconds. Preprocessing time is excluded; recorded times below $10^{-3}$ seconds are displayed at the left boundary. PS, DS, BR, and HB denote primal simplex, dual simplex, barrier, and homogeneous barrier, respectively. The first-order solvers \SCS{}, \ABIP{}, \PDLP{}, and \HPDLP{} are shown with solid lines.}
    \label{fig:infeas-detection}
\end{figure}

At $\varepsilon=10^{-8}$, \HPDLP{} detects 326 of the 334 instances in \MIPLIBInfeas{} without presolve and all 334 instances with \PSLP{}. It also detects all 29 infeasible \netlib{} instances in both presolve settings. \autoref{fig:infeas-detection} shows the fraction of detected instances as a function of solve time.

With a looser optimality tolerance, a run may satisfy \eqref{eq:optimality-tolerance} before a certificate is accepted and therefore terminate as Opt. On \MIPLIBInfeas{}, this occurs for 41 instances without presolve and 29 with \PSLP{} at $\varepsilon=10^{-4}$. Tightening $\varepsilon$ to $10^{-8}$ reduces these counts to 3 and 0, respectively, while the certificate tolerances remain fixed. On the infeasible \netlib{} collection, one instance terminates as Opt.\ in each presolve setting at $\varepsilon=10^{-4}$, compared with none at $\varepsilon=10^{-8}$.

These Opt.\ outcomes do not count as detections, but their actual termination times are included in SGM10. Therefore, the lower SGM10 values at $\varepsilon=10^{-4}$ should not be interpreted as faster infeasibility or unboundedness detection. Across both collections and both optimality tolerances, no \HPDLP{} run is classified as Fail.

\subsection{Feasible LP benchmarks}
\label{subsec:feas}

We also evaluate \HPDLP{} on the feasible LP benchmarks to assess its ability to solve LPs to the prescribed optimality tolerances. \autoref{tab:mip-relaxations-results} and \autoref{tab:mittelmann-results} report results alongside the CPU baselines and the GPU solvers \cuPDLPx{} and \HPRLPC{}.

\begin{table}[!htb]
\centering
\caption{Results on the 383 \MIPLIB{} relaxations. Outcome categories and SGM10 are defined in \autoref{subsec:term}. The GPU solvers at each optimality tolerance are compared separately with the CPU baselines. Highlighting follows \autoref{tab:miplib-infeas-results}, applied to the highest Solved count and lowest SGM10 within each presolve setting.}
\label{tab:mip-relaxations-results}
\footnotesize
\setlength{\tabcolsep}{4pt}
\begin{adjustbox}{width=0.8\linewidth}
\begin{threeparttable}
\begin{tabular}{lrrrrrrrr}
\toprule
& \multicolumn{4}{c}{Without presolve}
& \multicolumn{4}{c}{With \PSLP{}} \\
\cmidrule(lr){2-5}\cmidrule(lr){6-9}
Solver & Solved & Limit & Fail & SGM10 (sec) & Solved & Limit & Fail & SGM10 (sec) \\
\midrule
\Gurobi{} primal simplex & 358 & 25 & 0 & 28.82 & 366 & 17 & 0 & 22.96 \\
\Gurobi{} dual simplex   & 366 & 17 & 0 & 16.41 & 369 & 14 & 0 & 13.01 \\
\Gurobi{} barrier        & 361 & 1 & 21 & 13.82 & 368 & 1 & 14 & 8.77 \\
\Gurobi{} homogeneous barrier & 365 & 1 & 17 & 13.95 & 369 & 1 & 13 & 9.55 \\
\HiGHS{}                 & 353 & 30 & 0 & 39.03 & 360 & 23 & 0 & 26.30 \\
\MOSEK{}                 & \textbf{\underline{382}} & 1 & 0 & 10.07 & \textbf{\underline{382}} & 1 & 0 & 7.92 \\
\SCS{}                   & 308 & 57 & 18\tnote{a} & 149.13 & 323 & 46 & 14\tnote{b} & 114.51 \\
\ABIP{} ($10^{-4}$)      & 281 & 91 & 11\tnote{c} & 245.25 & 313 & 67 & 3 & 147.19 \\
\ORTools{} \PDLP{}       & 363 & 18 & 2\tnote{d} & 36.36 & 373 & 10 & 0 & 24.18 \\
\midrule
\HPRLPC{} ($10^{-4}$)    & 379 & 4 & 0 & 4.34 & 382 & 1 & 0 & 2.34 \\
\cuPDLPx{} ($10^{-4}$)   & 382 & 1 & 0 & 2.27 & 382 & 1 & 0 & 1.43 \\
\HPDLP{} ($10^{-4}$)   & \textbf{383} & 0 & 0 & \textbf{1.91} & \textbf{383} & 0 & 0 & \textbf{1.07} \\
\midrule
\HPRLPC{} ($10^{-8}$)    & 371 & 12 & 0 & 10.01 & 380 & 3 & 0 & 5.91 \\
\cuPDLPx{} ($10^{-8}$)   & 376 & 7 & 0 & 6.29 & 381 & 2 & 0 & \textbf{\underline{3.60}} \\
\HPDLP{} ($10^{-8}$)   & 379 & 4 & 0 & \textbf{\underline{5.55}} & \textbf{\underline{382}} & 1 & 0 & 4.02 \\
\bottomrule
\end{tabular}
\begin{tablenotes}
\footnotesize
\item[a] Includes 18 false infeasibility reports.
\item[b] Includes 14 false infeasibility reports.
\item[c] Includes 2 false infeasibility reports and 9 crashes.
\item[d] Includes 2 false infeasibility reports.
\end{tablenotes}
\end{threeparttable}
\end{adjustbox}
\end{table}

\begin{table}[!htb]
\centering
\caption{Results on the 49 instances of Mittelmann's LP benchmark. Column definitions and highlighting conventions follow \autoref{tab:mip-relaxations-results}.}
\label{tab:mittelmann-results}
\footnotesize
\setlength{\tabcolsep}{4pt}
\begin{adjustbox}{width=0.8\linewidth}
\begin{threeparttable}
\begin{tabular}{lrrrrrrrr}
\toprule
& \multicolumn{4}{c}{Without presolve}
& \multicolumn{4}{c}{With \PSLP{}} \\
\cmidrule(lr){2-5}\cmidrule(lr){6-9}
Solver & Solved & Limit & Fail & SGM10 (sec) & Solved & Limit & Fail & SGM10 (sec) \\
\midrule
\Gurobi{} primal simplex & 35 & 14 & 0 & 1434.07 & 40 & 9 & 0 & 956.45 \\
\Gurobi{} dual simplex   & 39 & 10 & 0 & 747.41 & 42 & 7 & 0 & 484.65 \\
\Gurobi{} barrier        & 42 & 0 & 7 & 90.72 & 45 & 0 & 4 & 50.38 \\
\Gurobi{} homogeneous barrier & 40 & 1 & 8 & 120.32 & 44 & 0 & 5 & 56.41 \\
\HiGHS{}                 & 34 & 15 & 0 & 1525.93 & 34 & 13 & 2 & 1548.74 \\
\MOSEK{}                 & 45 & 1 & 3 & 202.95 & 43 & 3 & 3\tnote{a} & 171.47 \\
\SCS{}                   & 24 & 19 & 6\tnote{b} & 2657.89 & 29 & 15 & 5\tnote{c} & 2169.78 \\
\ABIP{} ($10^{-4}$)      & 30 & 16 & 3 & 1617.78 & 35 & 11 & 3 & 1243.28 \\
\ORTools{} \PDLP{}       & 40 & 9 & 0 & 334.09 & 43 & 4 & 2\tnote{d} & 275.56 \\
\midrule
\HPRLPC{} ($10^{-4}$)    & \textbf{49} & 0 & 0 & 13.45 & \textbf{49} & 0 & 0 & 7.84 \\
\cuPDLPx{} ($10^{-4}$)   & 48 & 1 & 0 & 11.89 & \textbf{49} & 0 & 0 & 6.98 \\
\HPDLP{} ($10^{-4}$)   & \textbf{49} & 0 & 0 & \textbf{10.71} & \textbf{49} & 0 & 0 & \textbf{6.68} \\
\midrule
\HPRLPC{} ($10^{-8}$)    & \textbf{\underline{47}} & 2 & 0 & 43.03 & \textbf{\underline{48}} & 1 & 0 & 22.32 \\
\cuPDLPx{} ($10^{-8}$)   & 45 & 4 & 0 & \textbf{\underline{38.47}} & \textbf{\underline{48}} & 1 & 0 & \textbf{\underline{19.69}} \\
\HPDLP{} ($10^{-8}$)   & 45 & 4 & 0 & 38.97 & 47 & 2 & 0 & 22.29 \\
\bottomrule
\end{tabular}
\begin{tablenotes}
\footnotesize
\item[a] Includes 1 false infeasibility report and 2 out-of-memory terminations.
\item[b] Includes 5 false infeasibility reports and 1 crash.
\item[c] Includes 4 false infeasibility reports and 1 crash.
\item[d] Includes 2 false infeasibility reports.
\end{tablenotes}
\end{threeparttable}
\end{adjustbox}
\end{table}

At $\varepsilon=10^{-4}$, \HPDLP{} solves all instances in both benchmark sets, with and without presolve. At $\varepsilon=10^{-8}$, it solves 379 of the 383 \MIPLIB{} relaxations without presolve and 382 with \PSLP{}; the corresponding counts on Mittelmann's benchmark are 45 and 47 out of 49. The tighter tolerance increases SGM10 in all four benchmark--presolve settings, and some instances reach the time limit. Presolve with \PSLP{} reduces SGM10 at both tolerances and increases the solved counts at $\varepsilon=10^{-8}$.

\subsection{Iteration counts and computational resources}
\label{subsec:box-cone-experiments}

\begin{figure}[!htb]
    \centering
    \includegraphics[width=0.96\linewidth]{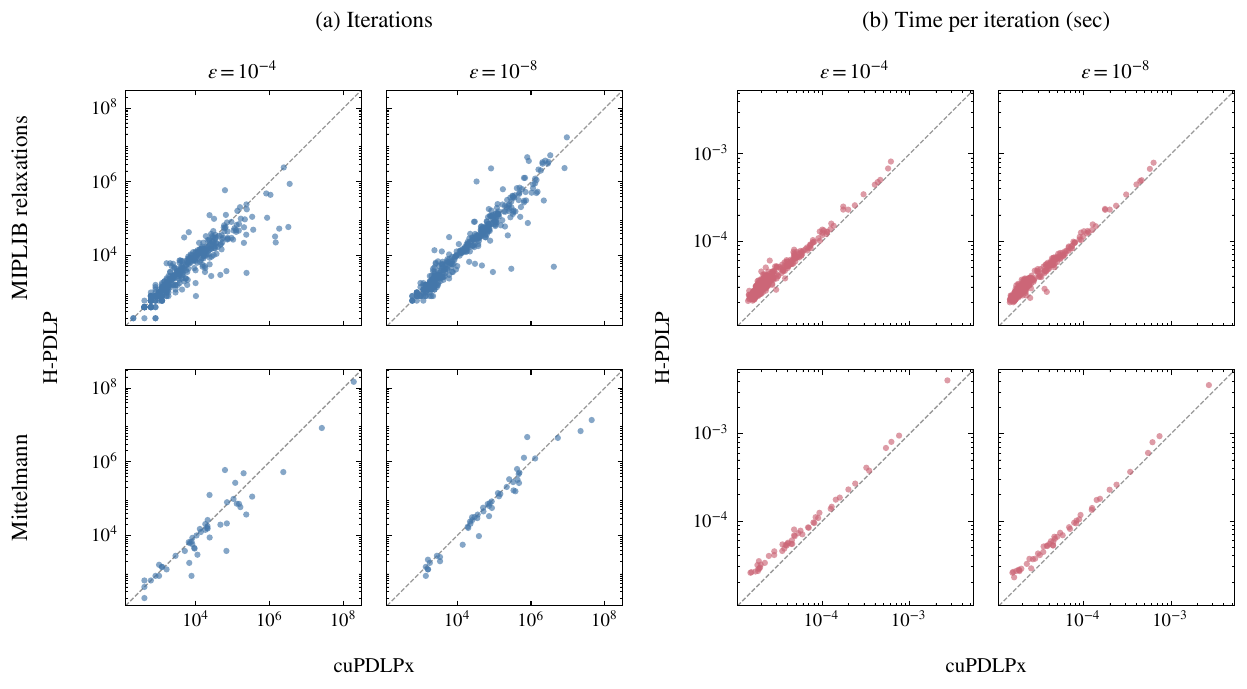}
    \caption{Iteration counts (a) and average time per iteration (b) for \HPDLP{} and \cuPDLPx{} with \PSLP{} preprocessing. Within each panel, rows correspond to \MIPLIB{} relaxations and Mittelmann, and columns to $\varepsilon=10^{-4}$ and $10^{-8}$. Only instances counted as Solved by both solvers and having positive iteration counts for both are included. Average time per iteration is computed as solve time divided by iteration count. Both axes are logarithmic; dashed lines indicate equality.}
    \label{fig:iteration-cost}
\end{figure}

To examine the computational costs underlying the preceding results, we use \cuPDLPx{} as a reference: both solvers use reflected Halpern iterations with adaptive restarting, and the implementation of \HPDLP{} builds on the \cuPDLPx{} codebase. Their different formulations and projection subproblems motivate a comparison of iteration counts and average time per iteration. \autoref{fig:iteration-cost} reports results with \PSLP{} preprocessing on feasible LP instances counted as Solved for both solvers under the outcome definitions in \autoref{subsec:term}. Across the two benchmark sets and two optimality tolerances, the geometric means of the per-instance \HPDLP{}/\cuPDLPx{} ratios range from 0.725 to 0.900 for iteration counts and from 1.308 to 1.447 for average time per iteration. Thus, lower aggregate iteration counts for \HPDLP{} are accompanied by higher average time per iteration.

\begin{figure}[!htb]
    \centering
    \includegraphics[width=0.9\linewidth]{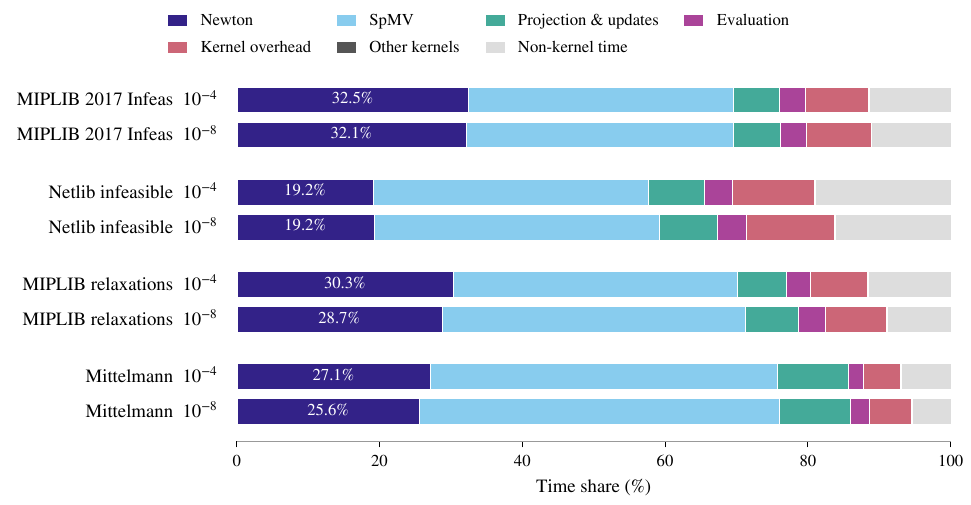}
   \caption{Computational time breakdown of \HPDLP{} with \PSLP{} preprocessing at $\varepsilon=10^{-4}$ and $10^{-8}$. Bars show mean per-instance time fractions over the recorded profiling intervals, with primal and dual contributions combined. ``Projection \& updates'' includes coordinate projection, reflection, and Halpern updates. ``Kernel overhead'' denotes the remaining time within kernels. ``Other kernels'' denotes other GPU kernels; non-kernel time is time spent outside CUDA kernels.}
    \label{fig:profile}
\end{figure}

We next examine the cost of the box-cone subproblems on all four benchmark sets with \PSLP{} preprocessing. Each PDHG step of \HPDLP{} solves primal and dual box-cone subproblems using a warm-started safeguarded Newton method. As shown in \autoref{tab:newton-iterations}, the medians of the per-instance mean Newton iteration counts range from 1.03 to 1.65, and the 90th percentiles range from 1.41 to 2.26. Despite these small typical counts, Newton solves account for approximately 19--33\% of the time within the recorded profiling intervals, averaged over instances, as shown in \autoref{fig:profile}. These solves introduce work beyond the coordinatewise projections used in \cuPDLPx{}, contributing to the higher average time per iteration of \HPDLP{}.

\begin{table}[!htb]
\centering
\caption{Newton iteration counts for the box-cone subproblems. Median and P90 are the median and 90th percentile of the per-instance mean counts, with equal weight per instance. Max is the largest recorded count in a single subproblem solve. Instances without recorded Newton evaluations are excluded.}
\label{tab:newton-iterations}
\footnotesize
\begin{adjustbox}{center}
\begin{tabular}{lcrrrrrr}
\toprule
\multicolumn{2}{c}{} & \multicolumn{3}{c}{Primal} & \multicolumn{3}{c}{Dual} \\
\cmidrule(lr){3-5}\cmidrule(lr){6-8}
Benchmark & $\varepsilon$ & Median & P90 & Max & Median & P90 & Max \\
\midrule
\MIPLIBInfeas{} & $10^{-4}$ & 1.45 & 2.06 & 48 & 1.65 & 2.26 & 11 \\
             & $10^{-8}$ & 1.44 & 2.07 & 48 & 1.65 & 2.26 & 11 \\
\netlib{} infeasible & $10^{-4}$ & 1.18 & 1.94 & 46 & 1.37 & 2.00 & 4 \\
                  & $10^{-8}$ & 1.16 & 1.89 & 46 & 1.35 & 2.00 & 4 \\
\midrule
\MIPLIB{} relaxations & $10^{-4}$ & 1.26 & 1.77 & 28 & 1.39 & 1.87 & 11 \\
                   & $10^{-8}$ & 1.09 & 1.41 & 72 & 1.15 & 1.47 & 16 \\
Mittelmann & $10^{-4}$ & 1.21 & 1.81 & 43 & 1.31 & 1.84 & 6 \\
           & $10^{-8}$ & 1.03 & 1.46 & 43 & 1.05 & 1.44 & 7 \\
\bottomrule
\end{tabular}
\end{adjustbox}
\end{table}

Beyond per-iteration time, \autoref{tab:host-memory} reports peak host-memory usage across the four benchmark sets. The reported peaks for \HPDLP{} remain below 14~GB on all four sets. For methods using sparse matrix factorizations, fill-in can substantially increase storage requirements. In these experiments, both \Gurobi{} barrier variants encounter out-of-memory terminations under the 128~GB limit on three benchmark sets, as does \MOSEK{} on Mittelmann; \SCS{}, using a direct linear system solver, reaches 102.14~GB on Mittelmann.

\begin{table}[!htb]
\centering
\caption{Peak host-memory usage (GB) across the four benchmark sets. Peaks are taken over runs with and without \PSLP{} preprocessing. Reported memory usage includes page cache and thus provides an upper bound on solver host-memory usage; GPU memory is excluded. Entries marked with an asterisk indicate the memory limit rather than a measured peak.}
\label{tab:host-memory}
\footnotesize
\setlength{\tabcolsep}{5pt}
\begin{adjustbox}{width=0.9\linewidth}
\begin{threeparttable}
\begin{tabular}{lrrrr}
\toprule
Solver
& \MIPLIBInfeas{}
& \netlib{} infeasible
& \MIPLIB{} relaxations
& Mittelmann \\
\midrule
\Gurobi{} primal simplex
& 15.64 & 0.05 & 74.27 & 28.99 \\
\Gurobi{} dual simplex
& 15.94 & 0.05 & 103.59 & 33.11 \\
\Gurobi{} barrier
& 128\tnote{*} & 0.19 & 128\tnote{*} & 128\tnote{*} \\
\Gurobi{} homogeneous barrier
& 128\tnote{*} & 0.23 & 128\tnote{*} & 128\tnote{*} \\
\HiGHS{}
& 30.94 & 0.05 & 85.57 & 50.10 \\
\MOSEK{}
& 79.08 & 0.27 & 115.16 & 128\tnote{*} \\
\SCS{}
& 61.92 & 0.07 & 36.94 & 102.14 \\
\ABIP{} ($10^{-4}$)
& 50.31 & 0.10 & 49.34 & 76.55 \\
\ORTools{} \PDLP{}
& 12.79 & 0.09 & 18.98 & 21.56 \\
\midrule
\HPRLPC{}\tnote{\ensuremath{\dagger}}
& --- & --- & 13.49 & 16.73 \\
\cuPDLPx{}\tnote{\ensuremath{\dagger}}
& --- & --- & 18.46 & 20.43 \\
\HPDLP{}\tnote{\ensuremath{\dagger}}
& 8.70 & 0.19 & 13.17 & 12.77 \\
\bottomrule
\end{tabular}
\begin{tablenotes}[flushleft]
\footnotesize
\item[*] At least one run encountered an out-of-memory termination under the 128~GB host-memory limit.
\item[\ensuremath{\dagger}] The peaks at $\varepsilon=10^{-4}$ and $10^{-8}$ are close; the larger value is reported.
\item[] Dashes indicate that no measurement is reported.
\end{tablenotes}
\end{threeparttable}
\end{adjustbox}
\end{table}

\section{Conclusion}
\label{sec:concl}

We developed \HPDLP{}, a GPU-based first-order LP solver that combines a generalized PDHG operator for a box-cone homogeneous self-dual formulation with restarted reflected Halpern iteration. The formulation encodes primal-dual optimal solutions and certificates of primal or dual infeasibility in a single conic system without introducing explicit bound slacks. The two box-cone projections reduce to scalar root-finding problems, so each PDHG iteration retains the sparse matrix-vector structure suited to GPUs and requires no matrix factorization. For fixed step sizes and the basic restart rule, we also derived a sufficient initialization under which the iterates remain uniformly bounded away from the trivial zero solution.

Numerical experiments on \MIPLIB{} relaxations and Mittelmann's LP benchmark demonstrate competitive performance against the tested CPU and GPU solvers, with remaining gaps primarily at the stricter accuracy setting. We also introduced \MIPLIBInfeas{}, a reproducible collection of 334 primal-infeasible and feasible but unbounded LPs derived from \MIPLIB{} 2017, extending the available benchmarks for evaluating infeasibility detection at scale. Experiments on this collection and the \netlib{} infeasible LP collection demonstrate effective certificate recovery by \HPDLP{}. These results support generalized PDHG applied to a box-cone HSD formulation as a practical approach to GPU-based first-order linear programming.

\section*{Acknowledgments}
We thank Jinwen Yang for the initial inspiration and many helpful discussions. We also thank Zedong Peng for helpful discussions on GPU computing and CUDA implementation. 

Haihao Lu is partially supported by AFOSR Grant No. FA9550-24-1-0051, ONR Grant No. N000142412735, and a Sloan Research Fellowship. Bo Tang is partially supported by ONR Grant No. N000142412735. We gratefully acknowledge NVIDIA Academic Grant Program for providing GPU computing resources.




\bibliographystyle{amsplain}
\bibliography{reference}

\end{document}